\documentclass[12pt]{article}

\usepackage{amsfonts,amsmath,amssymb,amsthm}
\usepackage{mathtools}
\DeclareMathOperator{\sign}{sgn}
\IfFileExists{mathabx.sty}{\usepackage{mathabx}}{}

\usepackage{tikz}
\usetikzlibrary{arrows}

\def\CC{\mathbb C}

\newcommand{\N}{\mathbb{N}}
\newcommand{\Z}{\mathbb{Z}}

\newcommand{\R}{\mathbb{R}}
\newcommand{\C}{\mathbb{C}}
\newcommand{\eps}{\epsilon}

\newtheorem{theorem}{Theorem}
\numberwithin{theorem}{section} 

\newtheorem{lemma}[theorem]{Lemma}
\newtheorem{proposition}[theorem]{Proposition}
\newtheorem{corollary}[theorem]{Corollary}
\theoremstyle{definition}

\newtheorem{example}[theorem]{Example}

\newtheorem{question}[theorem]{Question}
\newtheorem{remark}[theorem]{Remark}
\numberwithin{equation}{section}

\begin{document}
\title{Sharp Operator Khintchine Inequalities\\and the $Z_2$ Property}
\author{Chian Yeong Chuah, Zhen-Chuan Liu and Tao Mei}
\date{}
\maketitle
\begin{abstract}
We study the $p=1$ operator Khintchine inequality associated with
orthonormal systems and establish two quantitative implications
relating its optimal constant $A_1$ to the $Z_2$ property.
For an orthonormal system $W$ with finite $Z_2(W)$, we prove
$A_1(W)\leq\sqrt{1+Z_2(W)}$.
Conversely, there exists an absolute constant $\delta>0$ such that,
for canonical group unitaries indexed by any subset $V$ of a discrete
group, $A_1(V)<\sqrt2+\delta$ implies $Z_2(V)\leq2$.
The converse is detected by $2\times2$ matrix coefficients supported
on at most six group elements.
For the classical Rademacher sequence, we determine the sharp operator
Khintchine constant $A_1=\sqrt2$, answering the question left open by
Haagerup and Musat~\cite{Haagerup2007}.
Since this sequence has $Z_2=2$, the conclusion $Z_2(V)\leq2$ is optimal.
\end{abstract}

\section*{Introduction}
The classical Khintchine inequality states that, for $0<p<\infty$,
the $L_p$ and $L_2$ \mbox{(quasi-)norms} of linear combinations of Rademacher
functions are equivalent. More precisely, there exist optimal constants
$\alpha_p,\beta_p>0$ such that
\begin{equation}\label{khint}\tag{Kp}
 \beta_p^{-1}\left(\int_0^1\left|\sum_{k=1}^n c_k\eps_k(t)\right|^pdt\right)^{1/p}
 \leq\left(\sum_{k=1}^n|c_k|^2\right)^{1/2}
 \leq\alpha_p\left(\int_0^1\left|\sum_{k=1}^n c_k\eps_k(t)\right|^pdt\right)^{1/p}
\end{equation}
for every $n\in\N$ and $c_1,\ldots,c_n\in\C$. Here
$\eps_k(t)=\sign(\sin(2^k\pi t))$, $k\geq1$, are the Rademacher
functions on $[0,1]$. Analogous inequalities hold for other
orthonormal systems, including lacunary characters
$e^{2\pi i2^kt}$ and independent standard Gaussian variables.
For other orthonormal systems $W$, we use $\alpha_p(W)$ and
$\beta_p(W)$ for the corresponding optimal complex scalar constants,
and omit $W$ when it is clear from the context.
We are interested in the relation between the algebraic and
combinatorial properties of an orthonormal system and the validity
of these inequalities.

For characters on the torus indexed by a set $V\subset\Z$, the
classical $Z_2$ condition is
\[
 Z_2(V)=\sup_{h\in\Z\setminus\{0\}}
       \#\{(a,b)\in V^2:b-a=h\}<\infty.
\]
This condition is closely related to the work of Zygmund and Rudin
on lacunary Fourier series; see \cite{Rudin1960,Zygmund1974}.
Expanding the fourth moment and applying Cauchy--Schwarz gives
$\beta_4(V)\leq(1+Z_2(V))^{1/4}$. H\"older's inequality then yields
\[
 \alpha_1(V)\leq\beta_4(V)^2\leq\sqrt{1+Z_2(V)}.
\]
For example, $Z_2(\{2^k:k\geq1\})=1$.

The operator-valued analogue has been studied extensively
\cite{LustPiquard1986,LustPiquard1991,Pisier2017,Cadilhac19,JQSZ24}.
Noncommutative Khintchine inequalities are fundamental tools in the
geometry of noncommutative $L_p$-spaces and operator spaces
\cite{Pisier2003,Pisier2009,Xu2008}, in noncommutative martingale theory
\cite{PisierXu1997}, and in the study of random matrices
\cite{Rudelson1999,Tropp2015}.
Their applications include complete embeddings into noncommutative
$L_1$-spaces \cite{Junge2006,Xu2008}, martingale square-function
estimates \cite{PisierXu1997}, and covariance estimation
\cite{Rudelson1999,Tropp2015}.
The $p=1$ endpoint is particularly closely connected with the geometry
of preduals of von Neumann algebras \cite{Junge2006,Xu2008}.

The formulation uses column and row coefficient spaces, whose definitions
are recalled in Section~\ref{sec:ColSpaces}.
\begin{theorem}[\cite{LustPiquard1986,LustPiquard1991,Pisier2017}]
\label{thm:nc-Khintchine}
Let $0<p<\infty$. There exist optimal constants $A_p,B_p>0$, depending only
on $p$, such that for every $d,n\in\N$ and
$C_1,\ldots,C_n\in M_d(\C)$,
\begin{equation}\label{eq:nc-Khintchine}\tag{OKp}
 B_p^{-1}\left(\int_0^1\operatorname{tr}
     \left|\sum_{k=1}^n C_k\eps_k(t)\right|^pdt\right)^{1/p}
 \leq\|(C_k)\|_{S_p^d(\ell_2)}
 \leq A_p\left(\int_0^1\operatorname{tr}
     \left|\sum_{k=1}^n C_k\eps_k(t)\right|^pdt\right)^{1/p}.
\end{equation}
\end{theorem}
The trace in this formula is the usual unnormalized matrix trace.
At $p=1$, the middle term is the sum norm of the column and row
trace-class spaces, rather than a Hilbert-space norm.
Lust-Piquard proved the operator inequality for $1<p<\infty$
\cite{LustPiquard1986}; Lust-Piquard and Pisier treated $p=1$
\cite{LustPiquard1991}; and Pisier and Ricard established the range
$0<p<1$ \cite{Pisier2017}.

For character systems and $p>2$, the scalar Khintchine property
$(K_p)$ is the classical $\Lambda(p)$ property. The inclusion
$\Lambda(p)\subset\Lambda(q)$ for $2<q<p$ follows from monotonicity
of the $L_p$ norms. Bourgain~\cite[Theorem~2]{Bourgain1989} proved
that these inclusions are strict: for every $p>2$, there is a
$\Lambda(p)$ subset of $\Z$ that fails to be a $\Lambda(r)$ set
for every $r>p$. Harcharras~\cite{Harcharras1999} introduced the
corresponding completely bounded property $\Lambda(p)_{\mathrm{cb}}$,
expressed for $p>2$ by the matrix-coefficient inequality $(OK_p)$.
Together with interpolation and duality, the extrapolation results
of Pisier and Pisier--Ricard~\cite{Pisier2009,Pisier2017} describe
the usual passage to lower positive exponents. The identity at $p=2$
is automatic by orthonormality and provides no extrapolation hypothesis.

We study two quantitative implications relating $A_1$ to the
$Z_2$ property. Here $\alpha_1(W)$ is the optimal complex scalar
constant introduced above, whereas $A_1(W)$ is the corresponding
optimal matrix-coefficient constant, uniformly over all matrix
dimensions. Both constants are taken over finitely supported
coefficient families on $W$.

For every infinite subset $V$ of a discrete abelian group, one has
$A_1(V)\geq\sqrt2$. This follows from the Haagerup--Musat
lower-bound argument \cite[Theorem~2.9 and Remark~2.10]{Haagerup2007},
whose extension to commuting unitaries is recalled in
Lemma~\ref{lem:commuting-lower}. Thus $\sqrt2$ is a universal
benchmark for infinite abelian group systems.

Using the strategy of Haagerup and Musat~\cite{Haagerup2007}, we prove in
Theorem~\ref{thm:Main} and Remark~\ref{L4version} that
\[
 A_1(W)\leq\sqrt{1+Z_2(W)},\qquad
 Z_2(W)=\max\{Z_2^c(W),Z_2^r(W)\},
\]
where $Z_2^c$ and $Z_2^r$ are the centered-product constants defined
in Section~\ref{sec:Z2}. For canonical group unitaries indexed by
$V\subset\Gamma$, these are the column and row difference counts.

In the reverse direction, we prove that there exists an absolute
constant $\delta>0$ such that, for every subset $V$ of a discrete group,
\[
 A_1(V)<\sqrt2+\delta\quad\Longrightarrow\quad Z_2(V)\leq2.
\]
Thus an operator Khintchine constant sufficiently close to the
benchmark $\sqrt2$ forces every nonidentity difference to have at
most two representations, in both the column and row senses.
An explicit admissible value of $\delta$ is provided by the threshold
in Theorem~\ref{thm:quantitative-converse}.

Corollary~\ref{cor:near-optimal-all-p} gives a quantitative improvement
in the opposite direction. For group systems, the strict condition
$A_1(V)<\sqrt2+\delta$ above implies $(OK_q)$ for every $0<q\leq4$,
with constants depending only on $q$. In particular, near-optimal
behavior at $p=1$ yields a completely bounded fourth-moment estimate.
This upward self-improvement uses the intermediate conclusion
$Z_2(V)\leq2$, rather than interpolation from $p=1$ alone.
Proposition~\ref{prop:all-p} establishes the general passage from
finite $Z_2$ to this range of exponents.

The proof of the converse detects every obstruction $Z_2(V)\geq3$
using at most six group elements and $2\times2$ matrix coefficients.
Theorem~\ref{thm:converse} first treats the threshold $\sqrt2$ with
the same finite tests; Theorem~\ref{thm:quantitative-converse}
quantifies the gap in those obstructions.

A further main result is the determination of the optimal operator
Khintchine constant for the classical Rademacher sequence:
\[
 A_1(\{\eps_k:k\geq1\})=\sqrt2;
\]
see Theorem~\ref{rad26:thm:main}. Haagerup and Musat obtained the
upper bound $\sqrt3$ and left the determination of the optimal constant
open \cite[Remark~2.15]{Haagerup2007}.
Since the Rademacher characters have $Z_2=2$, this result also shows
that the conclusion $Z_2(V)\leq2$ in the converse theorem cannot be
improved for arbitrary discrete groups.

Our proof combines Cadilhac's self-adjoint factorization
\cite[Theorem~4.2 and Remark~4.6]{Cadilhac19}
with an extremality argument and a matrix-valued adaptation of the
even-function estimate of Lata\l a--Oleszkiewicz
\cite{LatalaOleszkiewicz1994}. For an extremizing Rademacher sum,
a norming functional identifies its mean modulus with an expression
determined by the factorization. This identity allows the even-function
estimate to control the column--row coefficient norm without losing
the sharp constant. The sharp real Gaussian constant follows as well
from the central limit theorem; see
Corollary~\ref{cor:real-gaussian-sharp}.

Section~\ref{sec:examples} gives examples that distinguish the roles
of $Z_2$, scalar coefficients, and matrix coefficients. The classical
Rademacher sequence and the quantum Rademacher system built from
Pauli matrices on at least two tensor factors both have $Z_2=2$,
but their matrix constants satisfy $A_1=\sqrt2$ and
$A_1\geq\sqrt{8/3}$, respectively. The latter bound holds even for
self-adjoint matrix coefficients. Thus $Z_2$ alone does not determine
$A_1$. We also construct an infinite set $V\subset\Z$ with
\[
 \alpha_1(V)<1.262<\sqrt2,\qquad Z_2(V)=\infty;
\]
see Corollary~\ref{prop:scalar-unbounded-differences}.
This shows that the matrix-coefficient hypothesis in the converse
theorem cannot be replaced by its scalar counterpart.

For infinite subsets $V\subset\Z$, the universal lower bound
$A_1(V)\geq\sqrt2$ and Theorem~\ref{thm:Main} show that $Z_2(V)=1$
implies $A_1(V)=\sqrt2$. Our converse gives only $Z_2(V)\leq2$.
We therefore ask whether
\[
 A_1(V)=\sqrt2\quad\Longleftrightarrow\quad Z_2(V)=1
\]
for every infinite subset of the integers. Here $Z_2(V)=1$ is the
Sidon property in the additive sense of Erd\H{o}s. This question
and the necessity of the infinitude assumption are discussed in
Section~\ref{sec:integer-Sidon-question}.

\section{Preliminaries}
\subsection{Noncommutative probability spaces}
Let $(\mathcal M,\tau)$ be a von Neumann algebra equipped with a
normal faithful tracial state. For $0<p<\infty$, write
$L_p(\mathcal M,\tau)$ for the associated noncommutative $L_p$ space,
with (quasi-)norm
\[
 \|f\|_p=\bigl(\tau(|f|^p)\bigr)^{1/p}.
\]
We set $L_\infty(\mathcal M)=\mathcal M$ with the operator norm.
We omit the trace from the notation when it is clear from the context.
The same construction is available for a normal semifinite faithful
trace. In particular, the ordinary trace $\operatorname{tr}$ on
$B(H)$ gives the Schatten spaces $S_p(H)$; this trace is not generally
a state. We always equip the coefficient algebra $M_d(\C)$ with
$\operatorname{tr}$ and the probability algebra $\mathcal M$ with
$\tau$, so matrix-valued functions carry the trace
$\operatorname{tr}\otimes\tau$.

The Cauchy--Schwarz and H\"older inequalities extend to this setting.
In particular, for $f\in M_d(\C)\bar\otimes\mathcal M$,
\begin{equation}\label{CS0}
 |(\operatorname{id}_{M_d}\otimes\tau)(f)|^2
 \leq(\operatorname{id}_{M_d}\otimes\tau)(|f|^2),
\end{equation}
and
\begin{equation}\label{Holder}
 \|f\varphi\|_r\leq\|f\|_p\|\varphi\|_q,
 \qquad \frac1r=\frac1p+\frac1q,
\end{equation}
whenever the indicated norms are finite. Here $\operatorname{id}$
denotes an identity map and $I_d$ an identity matrix. When no confusion
is possible, we also write $\tau$ for the entry-wise map
$\operatorname{id}_{M_d}\otimes\tau$ and omit tensor symbols in
expressions such as $C_k\otimes x_k$. See \cite{Pisier2003} for a
survey of noncommutative $L_p$ spaces.

\subsection{Column and row spaces}\label{sec:ColSpaces}
For a finite sequence $(x_k)$ in $L_p(\mathcal M,\tau)$, define
\[
 \|(x_k)\|_{L_p(\mathcal M;\ell_2^c)}
 =\left\|\left(\sum_kx_k^*x_k\right)^{1/2}\right\|_p,
 \qquad
 \|(x_k)\|_{L_p(\mathcal M;\ell_2^r)}
 =\left\|\left(\sum_kx_kx_k^*\right)^{1/2}\right\|_p.
\]
For $0<p<\infty$, the column and row spaces are the completions of
the finite sequences under these (quasi-)norms. For $p=\infty$, we
use the same formulas with the operator norm; only finite sequences
will be needed here.
For $0<p<2$, the sum space
\[
 L_p(\mathcal M;\ell_2)
 =L_p(\mathcal M;\ell_2^c)+L_p(\mathcal M;\ell_2^r)
\]
is equipped with the sum norm (a quasi-norm when $p<1$)
\[
 \|(x_k)\|_{L_p(\mathcal M;\ell_2)}
 =\inf_{x_k=y_k+z_k}
   \bigl\{\|(y_k)\|_{L_p(\mathcal M;\ell_2^c)}
          +\|(z_k)\|_{L_p(\mathcal M;\ell_2^r)}\bigr\}.
\]
For $2\leq p\leq\infty$, use the intersection space and its norm
\[
 \|(x_k)\|_{L_p(\mathcal M;\ell_2)}
 =\max\bigl\{\|(x_k)\|_{L_p(\mathcal M;\ell_2^c)},
             \|(x_k)\|_{L_p(\mathcal M;\ell_2^r)}\bigr\}.
\]
When $\mathcal M=M_d(\C)$ with its ordinary trace, write $S_p^d$
and $S_p^d(\ell_2)$ for these spaces. We occasionally omit $d$ when
it is fixed. In particular,
\begin{align*}
 \|(C_k)\|_{S_\infty^d(\ell_2)}
 &=\max\left\{\left\|\sum_kC_k^*C_k\right\|_\infty^{1/2},
               \left\|\sum_kC_kC_k^*\right\|_\infty^{1/2}\right\},\\
 \|(C_k)\|_{S_1^d(\ell_2)}
 &=\inf_{\substack{C_k=X_k+Y_k\\X_k,Y_k\in M_d(\C)}}
   \left\{\operatorname{tr}\left(\sum_kX_k^*X_k\right)^{1/2}
          +\operatorname{tr}\left(\sum_kY_kY_k^*\right)^{1/2}\right\}.
\end{align*}
The dual of the sum space is the intersection of the corresponding
dual spaces. Thus, for a finite matrix family,
\begin{equation}\label{eq:coefficient-duality}
 \|(C_k)\|_{S_1^d(\ell_2)}
 =\max_{\|(D_k)\|_{S_\infty^d(\ell_2)}\leq1}
       \operatorname{Re}\sum_k\operatorname{tr}(C_kD_k^*).
\end{equation}
These standard identifications are discussed in \cite{Pisier2003}.

\subsection{Group von Neumann algebras}\label{sec:group-algebra}
Let $\Gamma$ be a discrete group with identity $e$. Its left regular
representation is given by
$\lambda_g\delta_h=\delta_{gh}$ on $\ell_2(\Gamma)$.
The group von Neumann algebra $\mathcal L(\Gamma)$ is the weak
operator closure of the linear span of the canonical unitaries
$\lambda_g$. It has the normal faithful tracial state
\[
 \tau(f)=\langle\delta_e,f\delta_e\rangle,
 \qquad \tau(\lambda_g)=\delta_{g,e}.
\]
Consequently these unitaries form an orthonormal family:
\begin{remark}
For $x,y\in\Gamma$,
$\tau(\lambda_x^*\lambda_y)=\tau(\lambda_{x^{-1}y})=\delta_{x,y}$.
Every $f\in L_2(\mathcal L(\Gamma))$ has a Fourier expansion
$f=\sum_{x\in\Gamma}c_x\lambda_x$ converging in $L_2$, with
$(c_x)\in\ell_2(\Gamma)$.
\end{remark}
When $\Gamma$ is abelian, the Fourier transform identifies these
spaces with the classical spaces on its compact dual group
$\widehat\Gamma$.
For a subset $V\subset\Gamma$, define
\begin{align}
 Z_2^c(V)&=\sup_{g\ne e}\#\{(x,y)\in V^2:x^{-1}y=g\},\notag\\
 Z_2^r(V)&=Z_2^c(V^{-1}),\qquad
 Z_2(V)=\max\{Z_2^c(V),Z_2^r(V)\}.\label{eq:group-Z2}
\end{align}
The pairs are ordered. In an abelian group the two counts agree.
We call $V$ a $Z_2$ set when $Z_2(V)<\infty$.

\subsection{Centered products for orthonormal systems}\label{sec:Z2}
Let $(\mathcal M,\tau)$ be a finite von Neumann algebra with a normal
faithful tracial state, and let $W=(x_i)_{i\in I}\subset\mathcal M$ be an
orthonormal system: $\tau(x_i^*x_j)=\delta_{ij}$.  Set
\[
 z^c_{ij}=x_i^*x_j-\delta_{ij}1,\qquad
 z^r_{ij}=x_ix_j^*-\delta_{ij}1.
\]
Define $Z_2^c(W)$ and $Z_2^r(W)$ as the least constants for which,
for every finitely supported scalar array $(a_{ij})$,
\begin{align}
 \left\|\sum_{i,j}a_{ij}z^c_{ij}\right\|_2^2
 &\leq Z_2^c(W)\sum_{i,j}|a_{ij}|^2,\label{eq:Z2-c}\tag{$Z_2^c$}\\
 \left\|\sum_{i,j}a_{ij}z^r_{ij}\right\|_2^2
 &\leq Z_2^r(W)\sum_{i,j}|a_{ij}|^2.\label{eq:Z2-r}\tag{$Z_2^r$}
\end{align}
Put $Z_2(W)=\max\{Z_2^c(W),Z_2^r(W)\}$.
We say that $W$ has the $Z_2$ property if $Z_2(W)<\infty$.
This definition extends the difference-multiplicity condition for
group unitaries and also makes sense for orthonormal systems in
$L_4(\mathcal M)$. The main theorem is first stated for bounded
systems; its $L_4$ extension is given in Remark~\ref{L4version}.

\begin{remark}\label{rk:discrete group}
For $W=\{\lambda_x:x\in V\}$ in a discrete group $\Gamma$, write
$z^c_{xy}=\lambda_x^*\lambda_y-\delta_{x,y}1$.
These centered products vanish when $x=y$ and equal
$\lambda_{x^{-1}y}$ otherwise. Grouping terms by their nonidentity
differences and using orthonormality gives, for every finitely
supported scalar array $(a_{xy})$,
\begin{equation}\label{eq:group-centered-parseval}
 \left\|\sum_{x,y\in V}a_{xy}z^c_{xy}\right\|_2^2
 =\sum_{h\in\Gamma\setminus\{e\}}
    \left|\sum_{\substack{x,y\in V\\x^{-1}y=h}}a_{xy}\right|^2.
\end{equation}
Cauchy--Schwarz on each inner sum gives the upper bound with constant
$Z_2^c(V)$. Choosing coefficients equal to one on any finite collection
of pairs with the same nonidentity difference gives the matching
lower bound. The row calculation is identical, with pairs grouped
according to $xy^{-1}$. Equivalently, the column and row Gram
matrices are block diagonal, with all-ones blocks indexed by the
respective differences. Consequently,
\[
 Z_2^c(W)=Z_2^c(V),\qquad Z_2^r(W)=Z_2^r(V).
\]
Thus $Z_2(W)=\max\{Z_2^c(W),Z_2^r(W)\}=Z_2(V)$, so the centered-product
formulation recovers the usual difference-multiplicity condition.
\end{remark}

\section{Near-optimal inequalities imply the $Z_2$ property}\label{sec:converse}
\subsection{The operator Khintchine inequality at $p=1$}

We use the coefficient norm $S_1(\ell_2)$ defined in
Section~\ref{sec:ColSpaces}, with the usual nonnormalized matrix trace
$\operatorname{tr}$. Lust-Piquard and Pisier proved the following
operator analogue of the classical Khintchine inequality
\cite{LustPiquard1991}; Haagerup and Musat determined its optimal
constant \cite{Haagerup2007}.
\begin{theorem}[\cite{LustPiquard1991,Haagerup2007}]\label{thm:ncKh}
Let $V=\{3^j:j\geq0\}$. For every $d\geq1$ and every finitely
supported family $(C_k)_{k\in V}\subset M_d(\mathbb C)$,
\begin{equation}\label{eq:ncKh}\tag{OK1}
\int_0^1\operatorname{tr}\left|\sum_{k\in V}C_ke^{2\pi ikt}\right|\,dt
\leq \|(C_k)\|_{S_1(\ell_2)}
\leq \sqrt2\int_0^1\operatorname{tr}\left|\sum_{k\in V}C_ke^{2\pi ikt}\right|\,dt.
\end{equation}
The constant $\sqrt2$ is optimal uniformly in $d$ and the coefficient family.
\end{theorem}
The corresponding inequality for $p>1$ was first proved by
Lust-Piquard \cite{LustPiquard1986}. The first inequality in
\eqref{eq:ncKh} holds for every subset of integers. As the next lemma
shows, any constant in the second inequality for an infinite subset of
integers must be at least $\sqrt2$.

\subsection{Zygmund's $Z_2$ property and the converse}

Recall the column and row difference counts $Z_2^c(V)$ and $Z_2^r(V)$
from Section~\ref{sec:group-algebra}, and write $Z_2(V)=\max\{Z_2^c(V),Z_2^r(V)\}$.
Theorem~\ref{thm:Main} will imply that the canonical group unitaries
indexed by $V$ satisfy the operator Khintchine inequality with
$A_1\leq\sqrt{1+Z_2(V)}$. In particular, $Z_2(V)=1$ gives
$A_1\leq\sqrt2$. For an infinite subset of an abelian group, the
following lemma shows that this upper bound is sharp. For subsets of
the integers, the condition $Z_2(V)=1$ is equivalent to the Sidon
property in the sense of Erd\"os: each integer has at most one
representation $k+j$ with $k,j\in V$ and $k\leq j$.

Haagerup and Musat proved sharpness of the constant $\sqrt2$ for
$V=\{2^k:k\geq1\}$; see \cite[Theorem~2.9 and Remark~2.10]{Haagerup2007}.
Their lower-bound argument extends to every infinite subset of a
discrete abelian group. The following lemma records this extension
for the more general setting of commuting unitaries.

\begin{lemma}\label{lem:commuting-lower}
Let $(u_k)_{k\geq1}$ be an infinite sequence of commuting unitaries
in a finite tracial von Neumann algebra $(\mathcal N,\nu)$, where
$\nu(1)=1$. Suppose that a constant $A>0$ satisfies
\begin{equation}\label{eq:testKhintchin}
\|(C_k)\|_{S_1(\ell_2)}
\leq A\left\|\sum_k C_k\otimes u_k\right\|_
{L_1(M_d\overline\otimes\mathcal N,\operatorname{tr}\otimes\nu)}
\end{equation}
for every $d\geq1$ and every finitely supported family
$(C_k)\subset M_d(\mathbb C)$. Then $A\geq\sqrt2$.
\end{lemma}
\begin{proof}
Fix $n\geq1$. The abelian von Neumann algebra generated by
$u_1,\ldots,u_n$ can be represented as $L_\infty(\Omega,\mu)$
for a probability space $(\Omega,\mu)$, with $\nu$ represented by
integration \cite[Chapter III]{Takesaki2002}. Thus each $u_k$ is
represented by a function with $|u_k(\omega)|=1$ for almost every
$\omega\in\Omega$.

Let $s_1,\ldots,s_n$ be independent Steinhaus variables. Multiplication
of the individual coefficients by unimodular scalars preserves
$\|(C_k)\|_{S_1(\ell_2)}$. Applying \eqref{eq:testKhintchin} to
$(s_kC_k)$ therefore gives
\begin{equation}\label{eq:AveKhintchine}
\|(C_k)\|_{S_1(\ell_2)}
\leq A\int_\Omega\operatorname{tr}
\left|\sum_{k=1}^ns_ku_k(\omega)C_k\right|\,d\mu(\omega).
\end{equation}
Average over the normalized Haar measure $m$ on $\mathbb T^n$.
Fubini's theorem and invariance of $m$ under coordinatewise rotations
yield
\begin{align*}
\|(C_k)\|_{S_1(\ell_2)}
&\leq A\int_{\mathbb T^n}\int_\Omega\operatorname{tr}
\left|\sum_{k=1}^ns_ku_k(\omega)C_k\right|\,d\mu(\omega)\,dm(s)\\
&=A\int_\Omega\int_{\mathbb T^n}\operatorname{tr}
\left|\sum_{k=1}^ns_ku_k(\omega)C_k\right|\,dm(s)\,d\mu(\omega)\\
&=A\int_{\mathbb T^n}\operatorname{tr}
\left|\sum_{k=1}^ns_kC_k\right|\,dm(s).
\end{align*}
Since $n$ and $d$ are arbitrary, the sharp operator Khintchine
constant for independent Steinhaus variables gives $A\geq\sqrt2$
\cite[Remark~2.10]{Haagerup2007}.
\end{proof}

In particular, every infinite subset $V$ of a discrete abelian group
satisfies
\begin{equation}\label{eq:infinite-abelian-lower}
 A_1(V)\geq\sqrt2.
\end{equation}
Indeed, its canonical group unitaries commute, so the lemma applies
to any countably infinite subfamily. This includes every infinite
subset of $\Z$; no $Z_2$ assumption is required.

We now prove the converse at the threshold $\sqrt2$. It is enough
to test $2\times2$ matrix coefficients on at most six distinct
frequencies.
\begin{theorem}\label{thm:converse}
Let $V$ be a subset of a discrete group $\Gamma$. Suppose that, for
every $1\leq r\leq6$, every pairwise distinct
$g_1,\ldots,g_r\in V$, and every
$C_1,\ldots,C_r\in M_2(\mathbb C)$,
\begin{equation}\label{eq:Khintchin}
\|(C_m)_{m=1}^r\|_{S_1(\ell_2)}
\leq\sqrt2\,(\operatorname{tr}\otimes\tau)
\left|\sum_{m=1}^rC_m\otimes\lambda_{g_m}\right|.
\end{equation}
Then $Z_2(V)\leq2$.
\end{theorem}

\begin{lemma}\label{lem:Cnorm}
Suppose a finite family $(C_j)\subset M_d(\mathbb C)$ satisfies
\[
\sum_jC_j^*C_j=\alpha I_d,
\qquad \sum_jC_jC_j^*=\alpha I_d
\]
for some $\alpha>0$. Then
$\|(C_j)\|_{S_1(\ell_2)}=d\sqrt\alpha$.
\end{lemma}
\begin{proof}
The column norm gives the upper bound $d\sqrt\alpha$. For the lower
bound, set $D_j=C_j/\sqrt\alpha$. Both square-function sums for
$(D_j)$ equal $I_d$, so $\|(D_j)\|_{S_\infty(\ell_2)}=1$.
Duality gives
\[
\|(C_j)\|_{S_1(\ell_2)}
\geq\sum_j\operatorname{tr}(C_jD_j^*)
=\frac1{\sqrt\alpha}\operatorname{tr}\left(\sum_jC_jC_j^*\right)
=d\sqrt\alpha.
\]
\end{proof}

\begin{proof}[Proof of Theorem~\ref{thm:converse}]
Suppose first that $Z_2^c(V)\geq3$. Choose $h\neq e$ such that the
directed graph with vertices in $V$ and arrows $x\longrightarrow xh$
whenever $x,xh\in V$ has at least three arrows. Every vertex has at
most one incoming and one outgoing arrow. If $h^2=e$, the arrows
occur in two-cycles, so the graph contains two disjoint two-cycles.
Otherwise three arrows contain either a three-cycle or one of the
path configurations $1+1+1$, $2+1$, and $3$ in
Figure~\ref{fig:disjoint_lines}. A cycle of length greater than three
contains a three-edge path. Thus these five configurations exhaust
the possibilities.

\begin{figure}[htbp]
\centering
\begin{tabular}{cl}
$1+1+1$ &
\begin{tikzpicture}[baseline=-.5ex,>=stealth,scale=.8]
\draw[->] (0,0) node[left]{$a$} -- (1,0) node[right]{$ah$};
\draw[->] (2.4,0) node[left]{$b$} -- (3.4,0) node[right]{$bh$};
\draw[->] (4.8,0) node[left]{$c$} -- (5.8,0) node[right]{$ch$};
\end{tikzpicture}\\[9pt]
$2+1$ &
\begin{tikzpicture}[baseline=-.5ex,>=stealth,scale=.8]
\draw[->] (0,0) node[left]{$a$} -- (1,0);
\node at (1.3,0) {$ah$};
\draw[->] (1.7,0) -- (2.7,0) node[right]{$ah^2$};
\draw[->] (4.8,0) node[left]{$c$} -- (5.8,0) node[right]{$ch$};
\end{tikzpicture}\\[9pt]
$3$ &
\begin{tikzpicture}[baseline=-.5ex,>=stealth,scale=.8]
\draw[->] (0,0) node[left]{$a$} -- (1,0);
\node at (1.3,0) {$ah$};
\draw[->] (1.7,0) -- (2.7,0);
\node at (3.1,0) {$ah^2$};
\draw[->] (3.6,0) -- (4.6,0) node[right]{$ah^3$};
\end{tikzpicture}\\[9pt]
$C_3$ &
\begin{tikzpicture}[baseline=-.5ex,>=stealth,scale=.8]
\node (a) at (0,0) {$a$};
\node (b) at (1.8,.55) {$ah$};
\node (c) at (3.6,0) {$ah^2$};
\draw[->] (a) -- (b); \draw[->] (b) -- (c);
\draw[->] (c) to[bend left=15] (a);
\end{tikzpicture}\\[9pt]
$C_2+C_2$ &
\begin{tikzpicture}[baseline=-.5ex,>=stealth,scale=.8]
\node (a) at (0,0) {$a$}; \node (b) at (1.8,0) {$ah$};
\node (c) at (4.2,0) {$c$}; \node (d) at (6,0) {$ch$};
\draw[->] (a) to[bend left=20] (b);
\draw[->] (b) to[bend left=20] (a);
\draw[->] (c) to[bend left=20] (d);
\draw[->] (d) to[bend left=20] (c);
\end{tikzpicture}
\end{tabular}
\caption{The five directed configurations arising from three repeated differences.}
\label{fig:disjoint_lines}\label{fig:disjoint_lines21}
\label{fig:disjoint_lines3}\label{fig:triangle_cycle}
\label{fig:two_2cycles_directed}
\end{figure}

We first record the trace calculation common to the three path
configurations. Let $v=\lambda_g$ with $g\neq e$, let $w$ be a group
unitary, and set
\[
B=\begin{pmatrix}(1+v)/2\\ w/\sqrt2\end{pmatrix},
\qquad R=\begin{pmatrix}1&\lambda_h\end{pmatrix}.
\]
Since $RR^*=2$, the tracial polar-decomposition identity for
rectangular operators gives
\begin{align*}
(\operatorname{tr}\otimes\tau)|BR|
&=\sqrt2\,(\operatorname{tr}\otimes\tau)(BB^*)^{1/2}\\
&=\sqrt2\,\tau\bigl((B^*B)^{1/2}\bigr)
=\sqrt2\,\tau\left(\left(1+\frac{v+v^*}{4}\right)^{1/2}\right)
<\sqrt2.
\end{align*}
Indeed, $Q=1+(v+v^*)/4$ is positive and $\tau(Q)=1$.
Cauchy--Schwarz gives $\tau(Q^{1/2})\leq1$, with equality only if
\begin{equation}\label{eq:CS-f111}
Q^{1/2}=1.
\end{equation}
This would imply $v+v^*=0$, which is impossible: the two group
unitaries are either distinct and linearly independent, or equal and
have sum $2v\neq0$.

\medskip
\noindent\textbf{The configuration $1+1+1$.}
Here $a,ah,b,bh,c,ch$ are six distinct group elements. Choose, in this
order, the coefficients
\[
C_0=\tfrac12e_{11},\quad C_1=\tfrac12e_{12},\quad
C_2=\tfrac1{\sqrt2}e_{21},\quad C_3=\tfrac1{\sqrt2}e_{22},\quad
C_4=\tfrac12e_{11},\quad C_5=\tfrac12e_{12}.
\]
Then
\[
\sum_{j=0}^5C_j^*C_j=\sum_{j=0}^5C_jC_j^*=I_2,
\qquad \|(C_j)_{j=0}^5\|_{S_1(\ell_2)}=2
\]
by Lemma~\ref{lem:Cnorm}. The corresponding sum factors as
\[
f_{111}=\begin{pmatrix}(\lambda_a+\lambda_c)/2\\
\lambda_b/\sqrt2\end{pmatrix}R
=\lambda_a\begin{pmatrix}(1+\lambda_{a^{-1}c})/2\\
\lambda_{a^{-1}b}/\sqrt2\end{pmatrix}R.
\]
The common trace calculation, with $v=\lambda_{a^{-1}c}$ and
$w=\lambda_{a^{-1}b}$, gives
$(\operatorname{tr}\otimes\tau)|f_{111}|<\sqrt2$,
contradicting \eqref{eq:Khintchin}.

\medskip
\noindent\textbf{The configuration $2+1$.}
The five elements $a,ah,ah^2,c,ch$ are distinct. Assign to them
\[
D_0=\tfrac12e_{11},\quad
D_1=\tfrac12e_{12}+\tfrac1{\sqrt2}e_{21},\quad
D_2=\tfrac1{\sqrt2}e_{22},\quad
D_3=\tfrac12e_{11},\quad D_4=\tfrac12e_{12}.
\]
Again both square-function sums equal $I_2$, and
$\|(D_j)_{j=0}^4\|_{S_1(\ell_2)}=2$. The corresponding sum is
\[
f_{21}=\begin{pmatrix}(\lambda_a+\lambda_c)/2\\
\lambda_{ah}/\sqrt2\end{pmatrix}R
=\lambda_a\begin{pmatrix}(1+\lambda_{a^{-1}c})/2\\
\lambda_h/\sqrt2\end{pmatrix}R.
\]
Taking $v=\lambda_{a^{-1}c}$ and $w=\lambda_h$ yields
$(\operatorname{tr}\otimes\tau)|f_{21}|<\sqrt2$,
again contradicting \eqref{eq:Khintchin}.

\medskip
\noindent\textbf{The configuration $3$.}
The four elements $a,ah,ah^2,ah^3$ are distinct. Assign to them
\[
E_0=\tfrac12e_{11},\quad
E_1=\tfrac12e_{12}+\tfrac1{\sqrt2}e_{21},\quad
E_2=\tfrac1{\sqrt2}e_{22}+\tfrac12e_{11},\quad
E_3=\tfrac12e_{12}.
\]
Both square-function sums equal $I_2$, so
$\|(E_j)_{j=0}^3\|_{S_1(\ell_2)}=2$. This time
\[
f_3=\begin{pmatrix}(\lambda_a+\lambda_{ah^2})/2\\
\lambda_{ah}/\sqrt2\end{pmatrix}R
=\lambda_a\begin{pmatrix}(1+\lambda_{h^2})/2\\
\lambda_h/\sqrt2\end{pmatrix}R.
\]
Since $h^2\neq e$, the common calculation applies with
$v=\lambda_{h^2}$ and $w=\lambda_h$. Thus
$(\operatorname{tr}\otimes\tau)|f_3|<\sqrt2$, another contradiction.

For the remaining two configurations we use scalar coefficients,
embedded in $M_2(\mathbb C)$ as multiples of $e_{11}$. This embedding
preserves both the $\ell_2$ coefficient norm and the $L_1$ norm of
the scalar sum.

\medskip
\noindent\textbf{The configuration $C_3$.}
Here $h$ has order three. Choose coefficient $1$ at each of
$a,ah,ah^2$, so the coefficient norm is $\sqrt3$. The operator
\[
P_3=\tfrac13(1+\lambda_h+\lambda_{h^2})
\]
is a projection of trace $1/3$. Consequently,
\[
f_{C_3}=\lambda_a+\lambda_{ah}+\lambda_{ah^2}=3\lambda_aP_3,
\qquad \tau(|f_{C_3}|)=3\tau(P_3)=1.
\]
This contradicts \eqref{eq:Khintchin}, since $\sqrt3>\sqrt2$.

\medskip
\noindent\textbf{The configuration $C_2+C_2$.}
Here $h^2=e$, $h\neq e$, and $c\notin\{a,ah\}$, so the two cycles
are disjoint. Choose coefficient $1$ at each of $a,ah,c,ch$.
The coefficient norm is $2$, and
\[
f_{C_2+C_2}=2\lambda_a(1+u)P_2,
\qquad u=\lambda_{a^{-1}c},\qquad P_2=\tfrac12(1+\lambda_h).
\]
The operator $P_2$ is a projection of trace $1/2$. Set
$A=(1+u)P_2$. Since $a^{-1}c\notin\{e,h\}$,
\[
\tau(A^*A)=\tau\bigl((2+u+u^*)P_2\bigr)=1.
\]
Moreover, $|A|P_2=|A|$. Cauchy--Schwarz therefore gives
\[
\tau(|f_{C_2+C_2}|)=2\tau(|A|P_2)
\leq2\tau(A^*A)^{1/2}\tau(P_2)^{1/2}=\sqrt2.
\]
Equality would require $|A|=\sqrt2P_2$. Squaring yields
$P_2(u+u^*)P_2=0$. This is impossible, since expansion in the
group basis gives a nonempty sum with positive coefficients, even
when some of its group elements coincide. Hence the last inequality
is strict, contradicting \eqref{eq:Khintchin}.

All five configurations have been excluded, so $Z_2^c(V)\leq2$.
Finally, the hypothesis \eqref{eq:Khintchin} passes from $V$ to
$V^{-1}$: replace every coefficient by its adjoint and take the
adjoint of the sum. Both the $S_1(\ell_2)$ coefficient norm and the
$L_1$ norm are unchanged. Applying the column conclusion to $V^{-1}$
gives $Z_2^r(V)=Z_2^c(V^{-1})\leq2$, and hence $Z_2(V)\leq2$.
\end{proof}

\begin{theorem}[Quantitative converse]
\label{thm:quantitative-converse}
Put
\[
 J=\frac1{2\pi}\int_0^{2\pi}\sqrt{1+\tfrac12\cos t}\,dt,
 \qquad c_*:=\frac{\sqrt2}{J}=1.4381696565\ldots .
\]
Let $V$ be a subset of a discrete group $\Gamma$.
If $A_1(V)<c_*$, then $Z_2(V)\leq2$.
More generally, it suffices that there be a constant $K<c_*$ such that
\[
 \|(C_m)_{m=1}^r\|_{S_1(\ell_2)}
 \leq K\left\|\sum_{m=1}^r C_m\otimes\lambda_{g_m}\right\|_1
\]
for every $1\leq r\leq6$, every pairwise distinct
$g_1,\ldots,g_r\in V$, and every $C_1,\ldots,C_r\in M_2(\C)$.
Here $A_1(V)$ is the least constant when all matrix dimensions and
all finitely supported coefficient families are allowed.
\end{theorem}

\begin{proof}
Assume that $Z_2^c(V)\geq3$ and use the five configurations in the
proof of Theorem~\ref{thm:converse}. For $z\in\mathbb T$, replace
$(C_4,C_5)$ by $(zC_4,zC_5)$ in configuration $1+1+1$, and
$(D_3,D_4)$ by $(zD_3,zD_4)$ in configuration $2+1$.
In the three-edge path, keep $E_0,E_1$ fixed and replace only
\[
 E_2(z)=\frac1{\sqrt2}e_{22}+\frac z2e_{11},
 \qquad E_3(z)=\frac z2e_{12}.
\]
Thus the $e_{22}$ part of $E_2$ is unchanged.
In all three cases both square-function sums remain $I_2$, and the
coefficient norm remains $2$.  After multiplication by a group
unitary, the resulting sum has the form
\[
 F_z=B_zR,\qquad
 B_z=\begin{pmatrix}(1+zv)/2\\ w/\sqrt2\end{pmatrix},\qquad
 R=\begin{pmatrix}1&\lambda_h\end{pmatrix},
\]
for group unitaries $v,w$.  Since $RR^*=2$, the rectangular polar
decomposition identity gives
\[
 \|F_z\|_1=\sqrt2\,\tau\!\left[
 \left(1+\frac{zv+\overline z v^*}{4}\right)^{1/2}\right].
\]
If $\rho$ is the spectral distribution of $v$, averaging over $z$ and
using Haar invariance yields
\[
 \int_{\mathbb T}\tau\!\left[
 \left(1+\frac{zv+\overline z v^*}{4}\right)^{1/2}\right]dm(z)
 =J.
\]
Thus some $z$ gives $\|F_z\|_1\leq\sqrt2J$, and the coefficient-to-norm ratio is at least $c_*$.  The three-cycle configuration gives
the larger ratio $\sqrt3>c_*$.  For two disjoint two-cycles, retain
the notation $P=(1+\lambda_h)/2$, $u=\lambda_{a^{-1}c}$ and
$A=(1+u)P$ from the preceding proof.  In the corner $PL(\Gamma)P$
with normalized trace $\nu=2\tau$, put $T=P(u+u^*)P$ and
$Q=A^*A/2=P+T/2$. Write $v=\lambda_h$, so that $v^2=1$.
Expanding the two copies of $P=(1+v)/2$ gives
\[
 T=\frac14\sum_{\substack{i,j\in\{0,1\}\\ \eta\in\{-1,1\}}}
             v^iu^\eta v^j.
\]
Before collecting equal group elements, this is a sum of eight
terms with coefficient $1/4$. All coefficients are nonnegative,
so coincidences can only increase their squared sum. Since $T=T^*$,
Parseval's identity gives
\[
 \tau(T^2)=\|T\|_2^2\geq8\left(\frac14\right)^2=\frac12,
 \qquad \nu(T^2)\geq1.
\]
We have $0\leq Q\leq2P$ and $\nu(Q)=1$. Functional calculus in the
corner, whose identity is $P$, gives
\[
 (Q^{1/2}-P)^2
 \geq\frac{(Q-P)^2}{(1+\sqrt2)^2}.
\]
Consequently,
\[
 1-\nu(Q^{1/2})
 =\frac12\nu\bigl((Q^{1/2}-P)^2\bigr)
 \geq\frac{\nu(T^2)}{8(1+\sqrt2)^2}
 \geq\frac1{8(1+\sqrt2)^2}.
\]
The scalar coefficient norm is $2$, and the norm of the corresponding
sum is $\sqrt2\,\nu(Q^{1/2})$. Its ratio is therefore at least
\[
 \frac{2}{\sqrt2\left(1-\frac1{8(1+\sqrt2)^2}\right)}
 =\frac{16}{4+5\sqrt2}=1.4452083\ldots>c_*.
\]
Each configuration contradicts the assumed local estimate. Applying
the same argument to $V^{-1}$ gives $Z_2^r(V)\leq2$, and hence
$Z_2(V)\leq2$.
\end{proof}

\begin{remark}
The number $c_*$ in Theorem~\ref{thm:quantitative-converse} is an
explicit threshold supplied by the preceding witnesses; we do not
assert that this threshold is optimal.
The conclusion $Z_2(V)\leq2$ in Theorems~\ref{thm:converse}
and~\ref{thm:quantitative-converse}, however, is optimal in general. Let $V=\{r_1,r_2\}$ consist of two independent Rademacher
variables, viewed as canonical unitaries of $\mathbb Z_2^2$.
Then $Z_2(V)=2$, since $r_1^{-1}r_2=r_2^{-1}r_1$. For trace-class
operators $A,B$,
\[
\|Ar_1+Br_2\|_{L_1}
=\tfrac12\bigl(\|A+B\|_{S_1}+\|A-B\|_{S_1}\bigr).
\]
On the other hand, the parallelogram identity and the trace norm of
a column block give
\begin{align*}
\|(A,B)\|_{S_1(\ell_2)}
&\leq\operatorname{tr}\bigl((A^*A+B^*B)^{1/2}\bigr)\\
&=\left\|\begin{pmatrix}
(A+B)/\sqrt2&0\\ (A-B)/\sqrt2&0
\end{pmatrix}\right\|_{S_1}\\
&\leq\frac1{\sqrt2}\bigl(\|A+B\|_{S_1}+\|A-B\|_{S_1}\bigr)
=\sqrt2\,\|Ar_1+Br_2\|_{L_1}.
\end{align*}
Thus this system satisfies \eqref{eq:Khintchin}. Section~\ref{rad26:sec:main}
extends the same bound to the infinite Rademacher sequence.
\end{remark}

\begin{corollary}[Near-optimal $L_1$ inequalities imply $(OK_p)$]
\label{cor:near-optimal-all-p}
Let $\Gamma$ be a discrete group and let $V\subset\Gamma$ satisfy
$A_1(V)<c_*$. For every $0<p\leq4$ there is a constant $K_p\geq1$,
depending only on $p$, such that
\begin{equation}\label{eq:near-optimal-all-p}
 K_p^{-1}\|(C_g)_{g\in V}\|_{S_p^d(\ell_2)}
 \leq\left\|\sum_{g\in V}C_g\otimes\lambda_g\right\|_p
 \leq K_p\|(C_g)_{g\in V}\|_{S_p^d(\ell_2)}
\end{equation}
for every $d\geq1$ and every finitely supported family
$(C_g)_{g\in V}\subset M_d(\C)$.
The $L_p$ (quasi-)norm is taken on
$M_d(\C)\overline\otimes\mathcal L(\Gamma)$ with respect to
$\operatorname{tr}\otimes\tau$.
The constants are independent of $\Gamma$, $V$, $d$, and the number
of nonzero coefficients. At $p=4$, one has the explicit estimates
\begin{equation}\label{eq:near-optimal-fourth}
 \|(C_g)\|_{S_4^d(\ell_2)}
 \leq\left\|\sum_{g\in V}C_g\otimes\lambda_g\right\|_4
 \leq3^{1/4}\|(C_g)\|_{S_4^d(\ell_2)}.
\end{equation}
\end{corollary}

\begin{proof}
Theorem~\ref{thm:quantitative-converse} gives $Z_2(V)\leq2$.
Apply Proposition~\ref{prop:all-p} below to the canonical group
unitaries indexed by $V$. Their $Z_2$ constant is at most $2$, so the
resulting constants depend only on $p$.
The explicit estimate at $p=4$ follows from
\eqref{eq:Z2-fourth-endpoint}.
\end{proof}

\section{The $Z_2$ property implies the Khintchine inequality}\label{sec:Z2-Khintchine}

We first prove the operator Khintchine inequality from the
centered-product estimates of Section~\ref{sec:Z2}. We then discuss
sharpness and deduce the inequalities for the full range $0<p\leq4$.

\begin{theorem}\label{thm:Main}
Let $(\mathcal M,\tau)$ be a finite von Neumann algebra with a normal
faithful tracial state. Let $W=(x_i)_{i\in I}\subseteq\mathcal M$ be an
orthonormal system for which the centered-product estimates
\eqref{eq:Z2-c}--\eqref{eq:Z2-r} hold with finite constants
$Z_2^c(W)$ and $Z_2^r(W)$. Put
$Z_2(W)=\max\{Z_2^c(W),Z_2^r(W)\}$.
For every $d,l\in\mathbb N$, every family of pairwise distinct elements
$x_1,\ldots,x_l\in W$, and every $C_{x_1},\ldots,C_{x_l}\in M_d(\mathbb C)$,
\begin{equation}
\begin{split}
 \left\|\sum_{i=1}^l C_{x_i}\otimes x_i\right\|_1
 &\leq \|(C_{x_i})_{i=1}^l\|_{S_1^d(\ell_2)}\\
 &\leq \sqrt{1+Z_2(W)}
       \left\|\sum_{i=1}^l C_{x_i}\otimes x_i\right\|_1,
\end{split}
\end{equation}
Here the $L_1$ norm on $M_d(\mathbb C)\otimes\mathcal M$ is taken with
respect to $\operatorname{tr}\otimes\tau$, with $\operatorname{tr}$ the usual nonnormalized matrix trace.
\end{theorem}

\subsection{Proof of Theorem~\ref{thm:Main}}

Throughout this subsection, we use the hypotheses of
Theorem~\ref{thm:Main} and fix pairwise distinct $x_1,\ldots,x_l\in W$.
Write
\[
 E_d=\mathrm{id}_{M_d}\otimes\tau,\qquad
 \widehat f(x)=E_d\bigl(f(I_d\otimes x^*)\bigr),
 \qquad f\in M_d(\mathbb C)\otimes\mathcal M.
\]
The map $E_d$ is positive and unital. We abbreviate $C_{x_i}$ to $C_i$
within the proof.

\begin{lemma}\label{lem:bdd}
For every $f\in M_d(\mathbb C)\otimes\mathcal M$,
\begin{equation}\label{eq:S1}
\begin{split}
 \sum_{i=1}^l\widehat f(x_i)^*\widehat f(x_i)
 &\leq E_d(f^*f)\leq\|f\|_\infty^2 I_d,\\
 \sum_{i=1}^l\widehat f(x_i)\widehat f(x_i)^*
 &\leq E_d(ff^*)\leq\|f\|_\infty^2 I_d.
\end{split}
\end{equation}
Consequently,
$\| (\widehat f(x_i))_{i=1}^l\|_{S_\infty^d(\ell_2)}\leq\|f\|_\infty$.
\end{lemma}

\begin{proof}
Set $p=\sum_{i=1}^l\widehat f(x_i)\otimes x_i$. Orthonormality and the
tracial property give
\[
 E_d\bigl((f-p)^*(f-p)\bigr)
 =E_d(f^*f)-\sum_{i=1}^l\widehat f(x_i)^*\widehat f(x_i)\geq0.
\]
Similarly,
\[
 E_d\bigl((f-p)(f-p)^*\bigr)
 =E_d(ff^*)-\sum_{i=1}^l\widehat f(x_i)\widehat f(x_i)^*\geq0.
\]
The remaining bounds follow from positivity and unitality of $E_d$.
\end{proof}

\begin{lemma}\label{lem:centered-fourth}
Let $f=\sum_{i=1}^l C_i\otimes x_i$, and put
$P=\sum_{i=1}^l C_i^*C_i$ and $Q=\sum_{i=1}^l C_iC_i^*$.
Then
\begin{align}
 E_d\bigl((f^*f)^2\bigr)
 &\leq P^2+Z_2^c(W)\sum_{i=1}^l C_i^*QC_i,\label{eq:fstarf}\\
 E_d\bigl((ff^*)^2\bigr)
 &\leq Q^2+Z_2^r(W)\sum_{i=1}^l C_iPC_i^* .\label{eq:ffstar}
\end{align}
\end{lemma}

\begin{proof}
Orthonormality gives $E_d(f^*f)=P$ and $E_d(ff^*)=Q$.  Put
\[
 H=f^*f-P\otimes1=\sum_{i,j}C_i^*C_j\otimes z^c_{ij}.
\]
To justify the matrix amplification of \eqref{eq:Z2-c}, index the
products $z^c_{ij}$ by $\alpha$ and set
$G_{\alpha\beta}=\tau((z^c_\alpha)^*z^c_\beta)$.
The scalar estimate says $0\leq G\leq Z_2^c(W)I$. Thus, for any
matrix family $(A_\alpha)$,
\[
 E_d\!\left(\left(\sum_\alpha A_\alpha\otimes z^c_\alpha\right)^*
                   \left(\sum_\beta A_\beta\otimes z^c_\beta\right)\right)
 =\sum_{\alpha,\beta}G_{\alpha\beta}A_\alpha^*A_\beta
 \leq Z_2^c(W)\sum_\alpha A_\alpha^*A_\alpha.
\]
Indeed, testing on a vector $\xi\in\C^d$ reduces this to the same
Gram estimate for the vectors $(A_\alpha\xi)_\alpha$.
Apply this with $A_{ij}=C_i^*C_j$. Since $H=H^*$, we obtain
\[
 E_d(H^2)\leq Z_2^c(W)\sum_{j=1}^l C_j^*QC_j.
\]
Since $E_d(H)=0$, this proves \eqref{eq:fstarf}.  Applying the same
argument to $f^*$ and the row products $z^r_{ij}$ proves
\eqref{eq:ffstar}.
\end{proof}

\begin{corollary}\label{eq:inequality}
If $\|(C_i)_{i=1}^l\|_{S_\infty^d(\ell_2)}\leq1$, then
\[
 E_d\bigl((f^*f)^2\bigr)\leq(1+Z_2^c(W))E_d(f^*f),
 \qquad
 E_d\bigl((ff^*)^2\bigr)\leq(1+Z_2^r(W))E_d(ff^*).
\]
\end{corollary}

\begin{proof}
The hypothesis gives $P,Q\leq I_d$ and
$\sum_iC_i^*QC_i\leq P$, $\sum_iC_iPC_i^*\leq Q$.
The two estimates follow from Lemma~\ref{lem:centered-fourth}.
\end{proof}

We now apply the truncation and iteration argument of
Haagerup--Musat~\cite{Haagerup2007} with the preceding fourth-moment bounds.
Set
\[
K=\sqrt{1+Z_2(W)}.
\]

\begin{lemma}\label{Norm}
If $\|(C_i)_{i=1}^l\|_{S_\infty^d(\ell_2)}\leq1$, there exists
$f\in M_d(\mathbb C)\otimes\mathcal M$ such that
\begin{equation}\label{eq:bdd}
 \|f\|_\infty\leq K/2
\end{equation}
and
\begin{equation}\label{eq:max}
 \|(C_i-\widehat f(x_i))_{i=1}^l\|_{S_\infty^d(\ell_2)}\leq1/2.
\end{equation}
\end{lemma}

\begin{proof}
Let $g=\sum_{i=1}^l C_i\otimes x_i$, and form the self-adjoint matrix
\[
 A=\begin{bmatrix}0&g^*\\g&0\end{bmatrix}.
\]
For $R>0$, define the continuous odd function
\[
 F_R(t)=\begin{cases}
 -R,&t\leq-R,\\
 t,&-R<t<R,\\
 R,&t\geq R.
 \end{cases}
\]
If $J=\operatorname{diag}(I_d\otimes1,-I_d\otimes1)$, then $JAJ=-A$.
Functional calculus and oddness give $JF_R(A)J=-F_R(A)$.
Since $F_R(A)$ is self-adjoint, it has the form
\[
 F_R(A)=\begin{bmatrix}0&f^*\\f&0\end{bmatrix}
\]
for some $f\in M_d(\mathbb C)\otimes\mathcal M$, with
\begin{equation}\label{eq:fleqr}
 \|f\|_\infty=\|F_R(A)\|\leq R.
\end{equation}
Put $G_R(t)=t-F_R(t)$. The scalar inequality
$|G_R(t)|\leq t^2/(4R)$ implies $G_R(A)^2\leq A^4/(16R^2)$.
Comparing the diagonal blocks yields
\begin{equation}\label{eq:FourierNorm}
\begin{split}
 (g-f)^*(g-f)&\leq (g^*g)^2/(16R^2),\\
 (g-f)(g-f)^*&\leq (gg^*)^2/(16R^2).
\end{split}
\end{equation}
Let $D_i=C_i-\widehat f(x_i)=\widehat{g-f}(x_i)$.
By Lemma~\ref{lem:bdd}, positivity of $E_d$, and
Corollary~\ref{eq:inequality},
\begin{align*}
 \sum_{i=1}^l D_i^*D_i
 &\leq E_d\bigl((g-f)^*(g-f)\bigr)
 \leq\frac{E_d((g^*g)^2)}{16R^2}\\
 &\leq\frac{1+Z_2^c(W)}{16R^2}\sum_{i=1}^l C_i^*C_i
 \leq\frac{1+Z_2^c(W)}{16R^2}I_d.
\end{align*}
The corresponding row estimate is
\[
 \sum_{i=1}^l D_iD_i^*
 \leq\frac{1+Z_2^r(W)}{16R^2}\sum_{i=1}^l C_iC_i^*
 \leq\frac{1+Z_2^r(W)}{16R^2}I_d.
\]
Taking $R=K/2$ proves \eqref{eq:bdd} and \eqref{eq:max}.
\end{proof}

By homogeneity, for every $S=(C_i)_{i=1}^l$ we may choose an element
$\phi(S)\in M_d(\mathbb C)\otimes\mathcal M$, with $\phi(0)=0$, such that
\begin{equation}\label{eq:Norm}
\begin{split}
 \|\phi(S)\|_\infty&\leq\frac K2\|S\|_{S_\infty^d(\ell_2)},\\
 \|S-(\widehat{\phi(S)}(x_i))_{i=1}^l\|_{S_\infty^d(\ell_2)}
 &\leq\frac12\|S\|_{S_\infty^d(\ell_2)}.
\end{split}
\end{equation}

\begin{lemma}\label{lem:bdd2}
For every family $(C_i)_{i=1}^l$, there exists
$h\in M_d(\mathbb C)\otimes\mathcal M$ such that
$\widehat h(x_i)=C_i$ for $1\leq i\leq l$ and
\[
 \|h\|_\infty\leq K\|(C_i)_{i=1}^l\|_{S_\infty^d(\ell_2)}.
\]
\end{lemma}

\begin{proof}
Set $S_0=(C_i)_{i=1}^l$ and, for $j\geq0$, define
\[
 S_{j+1}=S_j-(\widehat{\phi(S_j)}(x_i))_{i=1}^l.
\]
By \eqref{eq:Norm},
\begin{equation}\label{eq:converge}
\begin{split}
 \|S_j\|_{S_\infty^d(\ell_2)}
 &\leq2^{-j}\|S_0\|_{S_\infty^d(\ell_2)},\\
 \|\phi(S_j)\|_\infty
 &\leq K2^{-j-1}\|S_0\|_{S_\infty^d(\ell_2)}.
\end{split}
\end{equation}
Hence the series
\begin{equation}\label{eq:h1}
 h=\sum_{j=0}^{\infty}\phi(S_j)
\end{equation}
converges in operator norm and satisfies the required norm bound.
The telescoping identity
\[
 S_0=\sum_{j=0}^n(\widehat{\phi(S_j)}(x_i))_{i=1}^l+S_{n+1}
\]
and continuity of the coefficient map in Lemma~\ref{lem:bdd} show that
$(\widehat h(x_i))_{i=1}^l=S_0$.
\end{proof}

\begin{proof}[Completion of the proof of Theorem~\ref{thm:Main}]
Let $F=\sum_{i=1}^l C_i\otimes x_i$. For every
$h\in M_d(\mathbb C)\otimes\mathcal M$, the normal trace pairing gives
\begin{equation}\label{eq:normal-pairing}
 (\operatorname{tr}\otimes\tau)(Fh^*)
 =\sum_{i=1}^l tr\bigl(C_i\widehat h(x_i)^*\bigr).
\end{equation}
The duality between $S_1^d(\ell_2)$ and $S_\infty^d(\ell_2)$, together
with Lemma~\ref{lem:bdd}, therefore implies
\[
 \|F\|_1
 =\sup_{\|h\|_\infty\leq1}|(\operatorname{tr}\otimes\tau)(Fh^*)|
 \leq\|(C_i)_{i=1}^l\|_{S_1^d(\ell_2)}.
\]
Conversely, for any $(D_i)_{i=1}^l$ in the unit ball of
$S_\infty^d(\ell_2)$, Lemma~\ref{lem:bdd2} provides $h$ with
$\widehat h(x_i)=D_i$ and $\|h\|_\infty\leq K$. Thus
\[
 \left|\sum_{i=1}^l \operatorname{tr}(C_iD_i^*)\right|
 =|(\operatorname{tr}\otimes\tau)(Fh^*)|\leq K\|F\|_1.
\]
Taking the supremum over $(D_i)$ proves the other inequality.
\end{proof}

\begin{remark}[The $L_4$ version]\label{L4version}
Theorem~\ref{thm:Main} extends to orthonormal systems in
$L_4(\mathcal M)$ with finite $Z_2^c$ and $Z_2^r$.
Fix a finite subfamily $x_1,\ldots,x_l$ and choose bounded spectral
truncations $y_i^{(n)}\to x_i$ in $L_4$.
Their $L_2$ Gram matrices
$G_n=(\tau((y_i^{(n)})^*y_j^{(n)}))_{i,j}$ converge to $I_l$.
For all large $n$, set
\[
 x_j^{(n)}=\sum_{i=1}^l y_i^{(n)}(G_n^{-1/2})_{ij}.
\]
These bounded families are orthonormal and $x_j^{(n)}\to x_j$ in $L_4$.
H\"older's inequality gives convergence in $L_2$ of both families of
centered products. Their finite Gram matrices therefore converge in
operator norm, so their column and row $Z_2$ constants converge to
those of the chosen subfamily, which are at most $Z_2^c(W)$ and
$Z_2^r(W)$. Apply Theorem~\ref{thm:Main} to $(x_j^{(n)})$ and pass to
the limit; for each fixed matrix coefficient family, the sums converge
in $L_1$. Since the finite subfamily was arbitrary, this proves the
claimed extension with the same constants.
\end{remark}

\subsection{Sharpness and other exponents}

\begin{remark}[Sharp dependence on $Z_2$]\label{rem:sharp-Z2-dependence}
The function $\sqrt{1+Z_2(W)}$ in Theorem~\ref{thm:Main} cannot be
improved for the full class of bounded orthonormal systems.
Indeed, let $e$ be a projection in a probability algebra with
$\tau(e)=t\in(0,1]$, and put $x=t^{-1/2}e$ and $W=\{x\}$.
Then $\tau(x^*x)=1$ and
\[
 Z_2(W)=\left\|\frac et-1\right\|_2^2=\frac1t-1.
\]
For every matrix $C$, the one-term column--row sum norm is
$\|(C)\|_{S_1(\ell_2)}=\|C\|_1$, whereas
$\|C\otimes x\|_1=\sqrt t\,\|C\|_1$. Consequently,
\[
 A_1(W)=\frac1{\sqrt t}=\sqrt{1+Z_2(W)}.
\]
This example proves sharpness for general orthonormal systems, not
for the narrower class of group unitaries. It also makes no assertion
about the optimal threshold in Theorem~\ref{thm:quantitative-converse}.
\end{remark}

For completeness, we record the following consequence of the fourth-moment
estimate and the extrapolation results of Pisier and Pisier--Ricard.

\begin{proposition}[Khintchine inequalities from the $Z_2$ property]\label{prop:all-p}
Let $(\mathcal M,\tau)$ be a finite von Neumann algebra with a normal
faithful tracial state, and let $W\subset L_4(\mathcal M)$ be an
orthonormal system with $Z_2(W)<\infty$.
For every $0<p\leq4$ there is a constant $K_p$, depending only on
$p$ and $Z_2(W)$, such that, for every pairwise distinct
$x_1,\ldots,x_l\in W$ and every $C_1,\ldots,C_l\in M_d(\C)$,
\begin{equation}\label{eq:all-p-Khintchine}
 K_p^{-1}\|(C_i)\|_{S_p^d(\ell_2)}
 \leq\left\|\sum_i C_i\otimes x_i\right\|_p
 \leq K_p\|(C_i)\|_{S_p^d(\ell_2)}.
\end{equation}
At $p=4$, one has the explicit estimates
\begin{equation}\label{eq:Z2-fourth-endpoint}
 \|(C_i)\|_{S_4^d(\ell_2)}
 \leq\left\|\sum_i C_i\otimes x_i\right\|_4
 \leq(1+Z_2(W))^{1/4}\|(C_i)\|_{S_4^d(\ell_2)}.
\end{equation}
\end{proposition}

\begin{proof}
We first treat bounded systems. Put $Z=Z_2(W)$,
$F=\sum_i C_i\otimes x_i$, $P=\sum_i C_i^*C_i$, and
$Q=\sum_i C_iC_i^*$. Taking the trace in
Lemma~\ref{lem:centered-fourth} gives
\[
 \begin{aligned}
 \|F\|_4^4
 &\leq\operatorname{tr}(P^2)
       +Z_2^c(W)\operatorname{tr}(Q^2)\\
 &\leq(1+Z)\|(C_i)\|_{S_4^d(\ell_2)}^4.
 \end{aligned}
\]
Since $E_d(F^*F)=P$ and $E_d(FF^*)=Q$, contractivity of the
trace-preserving conditional expectation $E_d$ on $L_2$ gives
$\|P\|_2,\|Q\|_2\leq\|F\|_4^2$. This proves
\eqref{eq:Z2-fourth-endpoint}.

For clarity, interpolation for $2\leq p\leq4$ can be carried out
on ambient $L_p$ spaces. On an auxiliary Rademacher cube set
\[
 D_i(h)=\mathbb E_{\varepsilon}(\varepsilon_i h),\qquad
 \mathcal T h=\sum_iD_i(h)\otimes x_i.
\]
The Bessel estimates of Lemma~\ref{lem:bdd}, applied on this cube,
yield
\[
 \sum_iD_i(h)^*D_i(h)\leq\mathbb E_{\varepsilon}(h^*h),\qquad
 \sum_iD_i(h)D_i(h)^*\leq\mathbb E_{\varepsilon}(hh^*).
\]
Thus $\|\mathcal T\|_{L_2\to L_2}\leq1$ and, by the fourth-moment
bound just proved,
$\|\mathcal T\|_{L_4\to L_4}\leq(1+Z)^{1/4}$.
Complex interpolation of noncommutative $L_p$ spaces
\cite{Pisier2003} bounds $\mathcal T$ for $2\leq p\leq4$.
Apply it to $h=\sum_i\varepsilon_i C_i$ and then use
Theorem~\ref{thm:nc-Khintchine}. We obtain
\begin{equation}\label{eq:synthesis-two-four}
 \left\|\sum_i C_i\otimes x_i\right\|_p
 \leq b_p\|(C_i)\|_{S_p^d(\ell_2)}
 \qquad(2\leq p\leq4),
\end{equation}
where $b_p$ depends only on $p$ and $Z$.
The reverse inequality with constant one follows from the same
conditional-expectation argument, now on $L_{p/2}$.

If $4/3\leq p<2$, then $2<p'=p/(p-1)\leq4$.
Column--row duality and orthonormality give
\[
 \begin{aligned}
 \|(C_i)\|_{S_p^d(\ell_2)}
 &=\sup_{\|(D_i)\|_{S_{p'}^d(\ell_2)}\leq1}
      \left|\sum_i\operatorname{tr}(C_iD_i^*)\right|\\
 &=\sup_{\|(D_i)\|_{S_{p'}^d(\ell_2)}\leq1}
      \left|(\operatorname{tr}\otimes\tau)
         \left(F\left(\sum_iD_i\otimes x_i\right)^*\right)\right|\\
 &\leq b_{p'}\|F\|_p.
 \end{aligned}
\]
In particular, a lower Khintchine estimate holds at $p_0=4/3$.
Pisier's extrapolation theorem~\cite[Theorem~1.1]{Pisier2009}
then gives the lower estimate for $1\leq p<4/3$.
Its extension by Pisier and Ricard~\cite[Corollary~2.2]{Pisier2017}
gives the same conclusion for $0<p<1$.
These results apply to orthonormal systems in noncommutative $L_2$
and preserve independence of the matrix dimension and the number
of coefficients; their constants depend only on the exponents
and the starting bound.

For completeness, the upper estimate for $0<p\leq2$ needs no
$Z_2$ assumption. Operator concavity of $t\mapsto t^{p/2}$,
orthonormality, and $\tau(1)=1$ give the upper bound by either
the column norm or the row norm. Decomposing the coefficients
and using the triangle or quasi-triangle inequality yields
\[
 \|F\|_p\leq\chi_p\|(C_i)\|_{S_p^d(\ell_2)},\qquad
 \chi_p=\max\{1,2^{1/p-1}\};
\]
see also~\cite[Remark~1.4]{Pisier2017}.
This completes the bounded case. For systems in $L_4$, apply the
bounded orthonormal approximation in Remark~\ref{L4version}
to each finite subfamily and pass to the limit. The approximating
$Z_2$ constants stay bounded, and the sums converge in every
$L_p$ with $0<p\leq4$; the exact fourth-moment estimate also
passes to the limit.
\end{proof}

\section{The Rademacher sequence}\label{rad26:sec:main}

We return to the Rademacher functions $(\eps_k)$ from the introduction.
The associated orthonormal system is
$W=\{\eps_k:k\geq1\}$, for which
$W=W^*$, $Z_2^c(W)=Z_2^r(W)=2$, and hence $Z_2(W)=2$.
Haagerup and Musat proved that
$A_1(W)\leq\sqrt3$.  The next theorem improves this to the sharp
constant $\sqrt2$.
Throughout this section, the coefficient spaces and matrix trace are
those of Section~\ref{sec:ColSpaces}.

\begin{theorem}\label{rad26:thm:main}
Let $d,n\in\N$ and $C_1,\ldots,C_n\in M_d(\C)$. Then
\begin{equation}\label{rad26:eq:main}
 \int_0^1\operatorname{tr}\left|\sum_{k=1}^n C_k\eps_k(t)\right|\,dt
 \leq \|(C_k)\|_{S_1^d(\ell_2)}
 \leq \sqrt{2}\int_0^1\operatorname{tr}
       \left|\sum_{k=1}^n C_k\eps_k(t)\right|\,dt.
\end{equation}
The factor $\sqrt{2}$ is optimal, so $A_1=\sqrt{2}$ in
\eqref{eq:nc-Khintchine} for $p=1$.
\end{theorem}

\subsection{Proof of Theorem~\ref{rad26:thm:main}}

We identify the probability algebra generated by
$\eps_1,\ldots,\eps_n$ with the finite cube
$\Omega_n=\{-1,1\}^n$, equipped with uniform probability measure,
so that $\eps_k(\omega)=\omega_k$. We use $\tau$ for integration on
this algebra and, as in the preliminaries, also for its entrywise
extension to matrix-valued functions. Thus
$\|f\|_{L_1}=\operatorname{tr}\tau(|f|)$, and this equals the
integral in \eqref{rad26:eq:main} when $f=\sum_k C_k\eps_k$.
The notation $\|X\|_2$ for a matrix refers to its Schatten $2$-norm
with the usual unnormalized trace.

For the column and row spaces introduced in
Section~\ref{sec:ColSpaces}, duality gives
\begin{equation}\label{rad26:eq:complex-dual}
 \|(C_k)\|_{S_1^d(\ell_2)}
 =\max_{\|(D_k)\|_{S_\infty^d(\ell_2)}\leq1}
       \operatorname{Re}\sum_k\operatorname{tr}(C_kD_k^*)
 =\max_{\substack{\sum_kD_k^*D_k\leq I_d\\
                         \sum_kD_kD_k^*\leq I_d}}
       \operatorname{Re}\sum_k\operatorname{tr}(C_kD_k^*).
\end{equation}
If the $C_k$ are self-adjoint, replacing $D_k$ by
$H_k=(D_k+D_k^*)/2$ preserves the real pairing. It also preserves
feasibility, since
\[
 \sum_k H_k^2
 \leq\frac12\sum_k(D_k^*D_k+D_kD_k^*)\leq I_d.
\]
Hence, for self-adjoint $C_k$,
\begin{equation}\label{rad26:eq:hermitian-dual}
 \|(C_k)\|_{S_1^d(\ell_2)}
 =\max_{\substack{D_k=D_k^*\\\sum_kD_k^2\leq I_d}}
       \sum_k\operatorname{tr}(C_kD_k).
\end{equation}
The self-adjoint decomposition below is contained in Cadilhac's
factorization theorem~\cite[Theorem~4.2 and Remark~4.6]{Cadilhac19}.
A self-contained finite-dimensional proof of the formulation for a
prescribed dual maximizer is given in Appendix~\ref{app:cadilhac-kkt}.

\begin{lemma}[Cadilhac's factorization, optimizer form]\label{rad26:lem:kkt}
Suppose that $C_1,\ldots,C_n\in M_d(\C)$ are self-adjoint and that
$(D_k)$ attains the maximum in \eqref{rad26:eq:hermitian-dual}.
Put $S=\sum_kD_k^2$. There exists $P\in M_d(\C)$, $P\geq0$, such that
\begin{equation}\label{rad26:eq:kkt}
 C_k=PD_k+D_kP\quad(1\leq k\leq n),\qquad SP=PS=P.
\end{equation}
Moreover, $P\neq0$ whenever $(C_k)\neq0$.
\end{lemma}

The next lemma adapts the even-function argument of
Lata\l a--Oleszkiewicz~\cite{LatalaOleszkiewicz1994} to the
matrix-valued absolute value.

\begin{lemma}\label{rad26:lem:even-modulus}
Let $C_1,\ldots,C_n\in M_d(\C)$ be self-adjoint, and put
$f=\sum_kC_k\eps_k$. Then
\begin{equation}\label{rad26:eq:lowerM}
 \operatorname{tr}\big((\tau(|f|))^2\big)
 \geq\frac12\tau\big(\operatorname{tr}(f^2)\big)
 =\frac12\sum_k\operatorname{tr}(C_k^2).
\end{equation}
\end{lemma}

\begin{proof}
We first recall the following estimate for self-adjoint matrices:
\begin{equation}\label{rad26:eq:modulus}
 \||H|-|K|\|_2\leq\|H-K\|_2.
\end{equation}
To verify it, let $(\lambda_a,e_a)$ and $(\mu_b,v_b)$ be orthonormal
eigenpairs for $H$ and $K$. Spectral expansion gives
\[
 \||H|-|K|\|_2^2
 =\sum_{a,b}(|\lambda_a|-|\mu_b|)^2
                  |\langle e_a,v_b\rangle|^2
 \leq\sum_{a,b}(\lambda_a-\mu_b)^2
                  |\langle e_a,v_b\rangle|^2
 =\|H-K\|_2^2.
\]

Now set $Y(\omega)=|f(\omega)|$. Since
$f(-\omega)=-f(\omega)$, we have $Y(-\omega)=Y(\omega)$.
For $\omega\in\Omega_n$, let $\omega^{(k)}$ denote the point
obtained by changing its $k$-th sign, and define
\[
 \partial_kY(\omega)=\frac{Y(\omega)-Y(\omega^{(k)})}{2}.
\]
Write the Walsh expansion
\[
 Y=\sum_{J\subseteq\{1,\ldots,n\}}\widehat Y(J)\eps_J,
 \qquad \eps_J=\prod_{k\in J}\eps_k,\qquad
 \widehat Y(J)=\tau(Y\eps_J).
\]
The evenness of $Y$ implies that $\widehat Y(J)=0$ for odd $|J|$.
Parseval's identity in the real Hilbert space of self-adjoint
matrices therefore gives
\begin{align*}
 \tau(\|Y\|_2^2)-\|\tau(Y)\|_2^2
 &=\sum_{J\neq\varnothing}\|\widehat Y(J)\|_2^2\\
 &\leq\frac12\sum_J|J|\,\|\widehat Y(J)\|_2^2
 =\frac12\sum_k\tau(\|\partial_kY\|_2^2).
\end{align*}
By \eqref{rad26:eq:modulus} and
$f(\omega)-f(\omega^{(k)})=2\eps_k(\omega)C_k$, we also have
$\|\partial_kY(\omega)\|_2\leq\|C_k\|_2$.
Since $\tau(\|Y\|_2^2)=\tau(\operatorname{tr}(f^2))
=\sum_k\operatorname{tr}(C_k^2)$, the preceding inequality
yields
\[
 \tau(\operatorname{tr}(f^2))
 -\operatorname{tr}((\tau(|f|))^2)
 \leq\frac12\tau(\operatorname{tr}(f^2)),
\]
which proves \eqref{rad26:eq:lowerM}.
\end{proof}

We now combine the two lemmas. The quadratic trace estimate in
Lemma~\ref{rad26:lem:even-modulus} alone does not give the required
column--row sum-norm inequality. Choosing an extremizing coefficient
family and a norming functional expresses its mean modulus in terms
of the positive multiplier from Lemma~\ref{rad26:lem:kkt}.
This is the step that transfers the estimate to the coefficient norm
without a dimension-dependent loss.

\begin{proof}[Proof of Theorem~\ref{rad26:thm:main}]
The first inequality in \eqref{rad26:eq:main} is already contained
in Theorem~\ref{thm:Main}. To prove the second, first fix $n,d$
and restrict to self-adjoint coefficients. The function
$(B_k)\mapsto\|\sum_kB_k\eps_k\|_{L_1}$ is a norm on this
finite-dimensional real vector space: if the sum vanishes, then
$B_k=\tau(\eps_k\sum_jB_j\eps_j)=0$ for every $k$.
Compactness of its unit sphere gives a nonzero sequence $(C_k)$
attaining
\begin{equation}\label{rad26:eq:A}
 A=\max_{\substack{B_k=B_k^*\\(B_k)\neq0}}
       \frac{\|(B_k)\|_{S_1^d(\ell_2)}}
            {\|\sum_{1\leq k\leq n}B_k\eps_k\|_{L_1}}.
\end{equation}
Here $A$ depends on the fixed values of $n,d$. Choose $(D_k)$
attaining \eqref{rad26:eq:hermitian-dual} for $(C_k)$, and put
$f=\sum_{k=1}^n C_k\eps_k$. By the choices of $(C_k)$ and $(D_k)$,
\[
 \|f\|_{L_1}=\frac1A\|(C_k)\|_{S_1^d(\ell_2)}
 =\frac1A\sum_k\operatorname{tr}(D_kC_k).
\]
We first establish the matrix identity
\begin{equation}\label{identity}
 \tau(|f|)=\frac1A\sum_kD_kC_k.
\end{equation}

Indeed, for every self-adjoint sequence $(B_k)$,
\begin{equation}\label{rad26:eq:support}
 \left|\sum_k\operatorname{tr}(D_kB_k)\right|
 \leq\|(B_k)\|_{S_1^d(\ell_2)}
 \leq A\left\|\sum_kB_k\eps_k\right\|_{L_1},
\end{equation}
with equality, and a positive pairing, at $(B_k)=(C_k)$.
On the real subspace of self-adjoint Rademacher sums, define
\[
 \ell\left(\sum_kB_k\eps_k\right)
 =\frac1A\sum_k\operatorname{tr}(D_kB_k).
\]
By \eqref{rad26:eq:support}, this functional has norm at most one.
The real Hahn--Banach theorem extends it to all self-adjoint
matrix-valued functions on $\Omega_n$ with the $L_1$ norm.
By duality, there exists a self-adjoint
function $U$ with $\|U(\omega)\|_\infty\leq1$, so that
$\ell(g)=\operatorname{tr}\tau(Ug)$. Testing on each coefficient
gives
\begin{equation}\label{rad26:eq:moments}
 \tau(\eps_kU)=\frac{D_k}{A}.
\end{equation}
For $f=\sum_k C_k\eps_k$, we obtain
\[
 \operatorname{tr}\tau(Uf)
 =\frac1A\sum_k\operatorname{tr}(D_kC_k)
 =\frac1A\|(C_k)\|_{S_1^d(\ell_2)}=\|f\|_{L_1}.
\]
At each point of the finite cube,
$\operatorname{tr}(U(\omega)f(\omega))\leq\operatorname{tr}|f(\omega)|$.
The nonnegative deficits have mean zero, so equality holds at every
point. Since $U=U^*$ and $\|U\|_\infty\leq1$, this forces
\[
 Uf=fU=|f|.
\]
Indeed, choose an orthonormal eigenbasis
$f(\omega)e_j=\lambda_j e_j$. For $\lambda_j\neq0$, equality in the
trace pairing forces
$\langle e_j,U(\omega)e_j\rangle=\sign(\lambda_j)$.
The contraction property then gives
$U(\omega)e_j=\sign(\lambda_j)e_j$.
Thus $U$ agrees with the sign of $f$ on its nonzero spectral
subspaces; its behavior on $\ker f(\omega)$ does not affect either
product. No differentiability of the trace norm is required.
Using~\eqref{rad26:eq:moments}, we conclude that
\begin{equation}\label{iden}
 \tau(|f|)=\tau(Uf)=\sum_k\tau(\eps_kU)C_k
 =\frac1A\sum_k D_kC_k,
\end{equation}
which proves~\eqref{identity}.

Apply Lemma~\ref{rad26:lem:kkt} to $(C_k)$ and $(D_k)$, and let
$P\geq0$ be the resulting nonzero multiplier. Define
\[
 S=\sum_kD_k^2,\qquad \Phi(X)=\sum_kD_kXD_k,\qquad
 T=\operatorname{tr}(P^2)+\operatorname{tr}(P\Phi(P))>0.
\]
The positivity of $T$ follows from $P\neq0$ and $\Phi(P)\geq0$.
We will use
\begin{equation}\label{rad26:eq:schwarz}
 \Phi(P)^2\leq\Phi(P^2),\qquad
 \operatorname{tr}(\Phi(P)^2)
 \leq\operatorname{tr}(\Phi(P^2))=\operatorname{tr}(P^2).
\end{equation}
To verify this, define $V:\C^d\to(\C^d)^n$ by
$V\xi=(D_1\xi,\ldots,D_n\xi)$. Then
$V^*V=S\leq I_d$, so $VV^*\leq I_{nd}$, and
\[
 \begin{aligned}
 \Phi(P)^2
 &=V^*(I_n\otimes P)VV^*(I_n\otimes P)V\\
 &\leq V^*(I_n\otimes P^2)V=\Phi(P^2).
 \end{aligned}
\]
Moreover, $SP=PS=P$ gives
$\operatorname{tr}(\Phi(P^2))=\operatorname{tr}(P^2S)
=\operatorname{tr}(P^2)$.

By \eqref{rad26:eq:kkt} and \eqref{iden}, we have
\begin{equation}\label{rad26:eq:mean-modulus}
 \begin{aligned}
 \tau(|f|)& =\frac1A\sum_kD_k(PD_k+D_kP)
   =\frac{\Phi(P)+SP}{A}
   =\frac{P+\Phi(P)}{A}.
 \end{aligned}
\end{equation}
Applying \eqref{rad26:eq:schwarz} implies
\[
 \operatorname{tr}((\tau(|f|))^2)
 =\frac{\operatorname{tr}(P^2)+2\operatorname{tr}(P\Phi(P))
              +\operatorname{tr}(\Phi(P)^2)}{A^2}
 \leq\frac{2T}{A^2}.
\]
On the other hand, orthogonality of the Rademacher functions and cyclicity of trace give
\begin{equation}\label{rad26:eq:energy}
 \begin{aligned}
 \tau(\operatorname{tr}(f^2))
 &=\sum_k\operatorname{tr}(C_k^2)
  =\sum_k\operatorname{tr}\big((PD_k+D_kP)^2\big)\\
 &=2\operatorname{tr}(P^2S)+2\operatorname{tr}(P\Phi(P))=2T.
 \end{aligned}
\end{equation}
Lemma~\ref{rad26:lem:even-modulus} now yields
\[
 \frac{2T}{A^2}\geq\operatorname{tr}((\tau(|f|))^2)\geq T.
\]

Since $T>0$, we obtain $A\leq\sqrt{2}$. This proves the desired
bound for all self-adjoint matrix coefficients and all $n,d$.

For arbitrary $C_k\in M_d(\C)$, consider the self-adjoint matrices
\[
 \widetilde C_k=\begin{pmatrix}0&C_k\\C_k^*&0\end{pmatrix}
 \in M_{2d}(\C).
\]
The trace norm of the Rademacher sum doubles, and we claim that
\begin{equation}\label{rad26:eq:dilation}
 \left\|\sum_k\widetilde C_k\eps_k\right\|_{L_1}
 =2\left\|\sum_k C_k\eps_k\right\|_{L_1},\qquad
 \|(\widetilde C_k)\|_{S_1^{2d}(\ell_2)}
 =2\|(C_k)\|_{S_1^d(\ell_2)}.
\end{equation}
To verify the second equality, for any decomposition $C_k=X_k+Y_k$
use
\[
 \widetilde C_k
 =\begin{pmatrix}0&X_k\\Y_k^*&0\end{pmatrix}
  +\begin{pmatrix}0&Y_k\\X_k^*&0\end{pmatrix}.
\]
The column norm of the first sequence and the row norm of the
second sequence both equal
$\|(X_k)\|_{S_1^d(\ell_2^c)}+\|(Y_k)\|_{S_1^d(\ell_2^r)}$.
Taking the infimum gives the inequality ``$\leq$''.
Conversely, if $(D_k)$ is feasible in
\eqref{rad26:eq:complex-dual}, then
$\widetilde D_k=\left(\begin{smallmatrix}0&D_k\\D_k^*&0\end{smallmatrix}\right)$
satisfies
\[
 \sum_k\widetilde D_k^2
 =\operatorname{diag}\left(\sum_kD_kD_k^*,\sum_kD_k^*D_k\right)
 \leq I_{2d}.
\]
Its pairing with $(\widetilde C_k)$ is
$2\operatorname{Re}\sum_k\operatorname{tr}(C_kD_k^*)$.
The dual formula therefore proves the reverse inequality in
\eqref{rad26:eq:dilation}. Applying the self-adjoint case in
dimension $2d$ and dividing by two proves \eqref{rad26:eq:main}.

Finally, take $d=1$, $n=2$ and $C_1=C_2=1$. Then
\[
 \|(C_1,C_2)\|_{S_1^1(\ell_2)}=\sqrt{2},\qquad
 \int_0^1|\eps_1(t)+\eps_2(t)|\,dt=1.
\]
Thus the constant $\sqrt{2}$ cannot be decreased.
\end{proof}

Identifying the Rademacher characters with the canonical generators
of $G=\bigoplus_{k\geq1}\mathbb Z/2\mathbb Z$, we obtain a family of group unitaries indexed by a subset $W\subset G$  with $Z_2(W)=2$ and optimal constant $A_1(W)=\sqrt{2}$.
Thus the conclusion $Z_2(W)\leq2$ in Theorem~\ref{thm:converse}
cannot in general be strengthened to $Z_2(W)\leq1$.

\section{Examples}\label{sec:examples}

\subsection{Group systems and classical Rademachers}

  \begin{example} This example shows that $Z_2^c$ and $Z_2^r$ can be different for nonabelian groups.
Let $(u_k)$ be free Haar unitaries, let $n\geq3$, and set
\[
 W_1=\{u_k:2\leq k\leq n\},\qquad
 W_2=\{u_1u_k:2\leq k\leq n\},\qquad W=W_1\cup W_2.
\]
Then $Z_2^c(W)=2$ and $Z_2^r(W)=n-1$.
Indeed, for distinct $i,j\geq2$,
$u_i^*u_j=(u_1u_i)^*(u_1u_j)$ gives the two column representations,
whereas $(u_1u_k)u_k^*=u_1$ has $n-1$ row representations.
Free reduction shows that no larger multiplicities occur.
For $n=2$, both second-order constants are $1$.
\end{example}
 \begin{example} (Discrete groups with the $Z_2$ property).
Let $\Gamma$ be a discrete group and $V\subset\Gamma$. With the
notation of Section~\ref{sec:group-algebra}, set $W=\lambda(V)$. The canonical group
unitaries form an orthonormal system.  For this system, the identities
in Remark~\ref{rk:discrete group} give $Z_2^c(W)=Z_2^c(V)$ and
$Z_2^r(W)=Z_2^r(V)$.  Hence $Z_2(W)=Z_2(V)$.
\begin{corollary}\label{thm:Main-Group}
Let $\Gamma$ be a discrete group and let $V\subset\Gamma$ satisfy
$Z_2(V)<\infty$. For pairwise distinct $x_1,\ldots,x_l\in V$ and
$C_{x_1},\ldots,C_{x_l}\in M_d(\C)$, put
$F=\sum_{i=1}^l C_{x_i}\otimes\lambda_{x_i}$. Then
\begin{equation}\label{eq:Szarek}
 \begin{aligned}
 \|F\|_{L_1(M_d\bar\otimes\mathcal L(\Gamma))}
   &\leq \|(C_{x_i})\|_{S_1^d(\ell_2)},\\
 \|(C_{x_i})\|_{S_1^d(\ell_2)}
   &\leq\sqrt{1+Z_2(V)}\,
      \|F\|_{L_1(M_d\bar\otimes\mathcal L(\Gamma))}.
 \end{aligned}
\end{equation}
Here the tensor product is equipped with the usual matrix trace and
the canonical tracial state on $\mathcal L(\Gamma)$.
\end{corollary}
\end{example}
\begin{example}[Length-lacunary subsets of free groups]\label{ex:length-lacunary}
Let $G$ be a free group equipped with the reduced-word length
associated with a free generating set. If $V=\{g_j:j\geq1\}$ satisfies
\[
 |g_{j+1}|>3|g_j|\qquad(j\geq1),
\]
then $Z_2^c(V)=Z_2^r(V)=1$.
\end{example}

\begin{proof}
First observe that $|h^2|\geq|h|$ for every $h\in G$.
Indeed, write $h=v w v^{-1}$ in cyclically reduced form, so that
$w$ is cyclically reduced and the displayed concatenation is reduced.
Then
\[
 |h^2|=2|v|+2|w|\geq2|v|+|w|=|h|.
\]
The identity element satisfies the same inequality.

Suppose that
\[
 a^{-1}b=c^{-1}d\neq e,\qquad a,b,c,d\in V.
\]
Let $t$ be the element of greatest length among the distinct elements
appearing in this equality, and put
\[
 L=|t|,\qquad
 M=\max\{|s|:s\in\{a,b,c,d\},\ s\neq t\}.
\]
Since each ordered pair has distinct entries, this maximum is
well-defined. The lacunarity hypothesis gives $L>3M$.
If $t$ occurs only once, the equality expresses $t$ as a product of
the other three entries or their inverses. Subadditivity then gives
$L\leq3M$, a contradiction.

Thus $t$ occurs twice, once in each ordered pair. If it occupies
the same position in both pairs, cancellation shows that
$(a,b)=(c,d)$. Otherwise, after exchanging the pairs if necessary,
the equality has the form
\[
 t^{-1}u=v^{-1}t,\qquad |u|,|v|\leq M.
\]
Set $h=v^{-1}t$. Then $t=vh$, $u=vh^2$, and $h^2=v^{-1}u$.
Consequently,
\[
 |t|\leq|v|+|h|
 \leq|v|+|h^2|
 \leq2|v|+|u|
 \leq3M,
\]
a contradiction. Hence every nonidentity column difference has
at most one representation. Applying the same argument to $V^{-1}$,
whose elements have the same lengths, proves the row assertion.
Since $V$ contains at least two elements, both constants equal one.
\end{proof}

\begin{corollary}\label{cor:length-lacunary}
Let $G$ and $(g_j)$ be as in Example~\ref{ex:length-lacunary}.
For $C_1,\ldots,C_l\in M_d(\C)$, put
$F=\sum_{j=1}^l C_j\otimes\lambda_{g_j}$. Then
\begin{equation}\label{eq:length-lacunary-Khintchine}
 \|F\|_{L_1(M_d\bar\otimes\mathcal L(G))}
 \leq\|(C_j)\|_{S_1^d(\ell_2)}
 \leq\sqrt2\,\|F\|_{L_1(M_d\bar\otimes\mathcal L(G))}.
\end{equation}
Here the tensor product carries the trace $\operatorname{tr}\otimes\tau$,
with the usual matrix trace and the canonical tracial state on
$\mathcal L(G)$.
\end{corollary}

\begin{proof}
Apply Corollary~\ref{thm:Main-Group} with $Z_2(V)=1$.
\end{proof}

\begin{remark}[Other length functions]\label{rem:length-lacunary-general}
The proof of Example~\ref{ex:length-lacunary} uses only symmetry,
subadditivity, and the inequality $|h^2|\geq|h|$.
Thus the same conclusion holds for a group equipped with a length
function $\ell:G\to[0,\infty)$ satisfying
\[
 \ell(e)=0,\quad \ell(h)>0\ (h\neq e),\quad
 \ell(h^{-1})=\ell(h),\quad \ell(hk)\leq\ell(h)+\ell(k),
\]
provided that $\ell(h^2)\geq\ell(h)$ for all $h\in G$ and
$\ell(g_{j+1})>3\ell(g_j)$ for all $j\geq1$.
In particular, the assertion applies to the usual length on $\Z$.
The square-length condition is an additional metric hypothesis,
not a consequence of torsion-freeness for arbitrary word metrics.
For example, in $\Z$ with generating set $\{\pm1,\pm2m\}$,
$m\geq2$, the associated length satisfies $\ell(m)=m$ but
$\ell(2m)=1$.

Strict doubling does not suffice even for the standard free-group
length. If $a$ is a free generator, the sequence
$a^{-1},a^3,a^7,a^{15},\ldots$ has successive lengths strictly
more than twice the preceding lengths, but
\[
 (a^{-1})^{-1}a^3=a^4=(a^3)^{-1}a^7.
\]
Its column difference constant is therefore at least two.
We do not claim that the factor three in
Example~\ref{ex:length-lacunary} is optimal.
\end{remark}

\begin{example}(Rademacher sequence)
    Let $W=\{\eps_k:k\geq1\}$ be the classical Rademacher system. Then
$W=W^*$ and
\[
 Z_2^c(W)=Z_2^r(W)=2,\qquad Z_2(W)=2.
\]
The sharp Khintchine constant is $A_1(W)=\sqrt2$; its proof is given
in Section~\ref{rad26:sec:main}.
\end{example} 
\subsection{Quantum Rademachers and Pauli subsystems}

\begin{example}(Quantum Rademacher sequences).
Let $\{\sigma_1, \sigma_2, \sigma_3\}$ denote the Pauli matrices, i.e.,
\[
\sigma_1 = \begin{pmatrix} 0 & 1 \\ 1 & 0 \end{pmatrix}, \quad
\sigma_2 = \begin{pmatrix} 0 & -i \\ i & 0 \end{pmatrix}, \quad
\sigma_3 = \begin{pmatrix} 1 & 0 \\ 0 & -1 \end{pmatrix}.
\]
For each $1 \leq k \leq N$ and $i \in \{1,2,3\}$, define
\[
\epsilon_{i,k} := I_2^{\otimes(k-1)} \otimes \sigma_i \otimes I_2^{\otimes(N-k)} \in M_{2^N}(\CC),
\]
the operator which acts as $\sigma_i$ on the $k$th qubit and as identity elsewhere.

Let
\[
\mathcal{A} = \{ \epsilon_{i,k} : 1 \leq k \leq N,\ 1 \leq i \leq 3 \} \subset M_{2^N}.
\]
With respect to the normalized trace $\tau(X) = 2^{-N} \operatorname{tr}(X)$,
the family $\mathcal{A}$ forms a noncommutative orthonormal system:
\[
\tau(\epsilon_{i,k}^* \epsilon_{j,l}) =
\begin{cases}
1 & \text{if } i = j, \ k = l, \\
0 & \text{otherwise.}
\end{cases}
\]
The system $\mathcal{A}$ is a quantum analogue of classical Rademacher variables.

A direct calculation gives
$Z_2^c(\mathcal{A})=Z_2^r(\mathcal{A})=2$.
Theorem~\ref{thm:Main} therefore gives the following estimate.
\begin{corollary}\label{cor:quantum-rademacher}
For any choice of $c_{i,k}\in M_d(\CC)$,
\[
\left\| (c_{i,k})_{i,k} \right\|_{S_1^d(\ell_2)}
\leq \sqrt{3}(\operatorname{tr}\otimes\tau)
\left(\left| \sum_{i,k} c_{i,k} \otimes \epsilon_{i,k} \right|\right),
\]
where $\operatorname{tr}$ is the usual unnormalized trace on $M_d(\CC)$
and $\tau=2^{-N}\operatorname{tr}$ is the normalized trace on
$M_{2^N}(\CC)$.
This is a noncommutative Khintchine inequality for quantum Rademacher functions.
\end{corollary}

The following example shows that, when $N\geq2$, the constant for this
system is strictly larger than $\sqrt{2}$, even for self-adjoint matrix coefficients.
For fixed $N$, let $A_1(\mathcal{A})$ denote the least constant such that
\[
 \|(c_{i,k})_{i,k}\|_{S_1^d(\ell_2)}
 \leq A_1(\mathcal{A})(\operatorname{tr}\otimes\tau)
 \left(\left|\sum_{i,k}c_{i,k}\otimes\epsilon_{i,k}\right|\right)
\]
holds for every $d\in\N$ and every family $c_{i,k}\in M_d(\CC)$.

\begin{proposition}\label{prop:quantum-rademacher-lower}
For the quantum Rademacher system $\mathcal{A}$ with $N\geq2$,
\begin{equation}\label{eq:quantum-rademacher-bounds}
 \sqrt{\frac{8}{3}}\leq A_1(\mathcal{A})\leq\sqrt{3}.
\end{equation}
The lower bound also holds if the coefficients are required to be
self-adjoint.
\end{proposition}

\begin{proof}
The upper bound follows from Corollary~\ref{cor:quantum-rademacher}.
For the lower bound, first take $N=2$ and use $M_4(\CC)$ as the
coefficient algebra. Set
\[
 c_{i,1}=\sigma_i\otimes I_2,\qquad
 c_{i,2}=-I_2\otimes\sigma_i,\qquad i=1,2,3.
\]
These six matrices are self-adjoint unitaries, so
\[
 \sum_{i=1}^3\sum_{k=1}^2c_{i,k}^*c_{i,k}
 =\sum_{i=1}^3\sum_{k=1}^2c_{i,k}c_{i,k}^*
 =6I_4.
\]
Lemma~\ref{lem:Cnorm} therefore gives
\begin{equation}\label{eq:quantum-rademacher-coefficient-norm}
 \|(c_{i,k})_{i,k}\|_{S_1^4(\ell_2)}=4\sqrt{6}.
\end{equation}
Consider
\[
 F=\sum_{i=1}^3\sum_{k=1}^2c_{i,k}\otimes\epsilon_{i,k},
 \qquad X=\sum_{i=1}^3\sigma_i\otimes\sigma_i.
\]
Regrouping the four tensor factors in $F$ in the order $(1,3,2,4)$
shows that $F$ is unitarily equivalent to
$X\otimes I_4-I_4\otimes X$.
Let $P$ be the flip on $\CC^2\otimes\CC^2$, defined by
$P(\xi\otimes\eta)=\eta\otimes\xi$. Direct calculation gives
\[
 X=2P-I_4.
\]
The symmetric and antisymmetric subspaces have dimensions $3$ and $1$,
respectively. Hence $X$ has eigenvalues $1$ with multiplicity $3$ and
$-3$ with multiplicity $1$. It follows that the eigenvalues of $F$ are
\[
 \begin{array}{c|ccc}
 \text{eigenvalue}&0&4&-4\\ \hline
 \text{multiplicity}&10&3&3
 \end{array}
\]
and thus $\operatorname{tr}_{16}(|F|)=24$, where
$\operatorname{tr}_{16}$ is the usual trace on $M_{16}(\CC)$.
Since the quantum trace on $M_4(\CC)$ is normalized,
\[
 (\operatorname{tr}\otimes\tau)(|F|)
 =\frac14\operatorname{tr}_{16}(|F|)=6.
\]
Combining this identity with
\eqref{eq:quantum-rademacher-coefficient-norm}, we obtain
\[
 A_1(\mathcal{A})\geq
 \frac{\|(c_{i,k})_{i,k}\|_{S_1^4(\ell_2)}}
 {(\operatorname{tr}\otimes\tau)(|F|)}
 =\frac{4\sqrt{6}}{6}=\sqrt{\frac{8}{3}}.
\]
For $N>2$, use the same coefficients at the first two sites and set
$c_{i,k}=0$ for $k>2$. The additional identity tensor factors do not
change the normalized quantum trace, so the same ratio is obtained.
All the coefficients used above are self-adjoint, which proves the
last assertion.
\end{proof}

\begin{theorem}[Pauli directional subsystems]
\label{thm:pauli-subsystems}
Let $\mathcal S\subseteq\{1,2,3\}\times\{1,\ldots,N\}$ be nonempty and
let $\mathcal A(\mathcal S)=\{\epsilon_{i,k}:(i,k)\in\mathcal S\}$.
If $|\mathcal S|=1$, then $A_1(\mathcal A(\mathcal S))=1$.  If
$|\mathcal S|\geq2$, then
\[
 A_1(\mathcal A(\mathcal S))=\sqrt2
 \quad\Longleftrightarrow\quad
 \text{at most one site }k\text{ contains all three directions.}
\]
If two sites are full, then
$A_1(\mathcal A(\mathcal S))\geq\sqrt{8/3}$.
\end{theorem}

\begin{proof}
Conjugation by the identity or a Pauli matrix at site $k$ changes the signs of the
three directions there by one of the four patterns
$(1,1,1),(1,-1,-1),(-1,1,-1),(-1,-1,1)$.
If at most one site is full, choose a global sign so that the desired
pattern at that site has product $1$; at each partial site, complete
the prescribed signs to a triple with product $1$.
Thus every sign pattern on the selected variables is implemented,
up to a global sign, by trace-preserving conjugation.

Write $q=2^N$ and enumerate the selected variables as $u_j\in M_q$.
Tensoring each coefficient with $u_j$ gives
\[
 \|(C_j\otimes u_j)\|_{S_1^{dq}(\ell_2)}
 =q\|(C_j)\|_{S_1^d(\ell_2)}.
\]
The upper bound follows by tensoring a column--row decomposition;
the reverse bound follows from \eqref{eq:coefficient-duality} by tensoring a
dual maximizing family with the same unitaries.
Apply the sharp Rademacher inequality of
Theorem~\ref{rad26:thm:main} to $(C_j\otimes u_j)$.
The sign symmetry makes every trace norm in its average equal to
$\operatorname{tr}_{dq}|\sum_j C_j\otimes u_j|$.
Dividing by $q$ proves $A_1\leq\sqrt2$.

For the lower bound, choose two distinct selected variables $u,v$ and
use the self-adjoint coefficients $C_1=u$, $C_2=v$ in $M_q$.
They satisfy $u^2=v^2=I_q$ and either commute or anticommute.
Consequently $X=u\otimes u$ and $Y=v\otimes v$ are commuting
self-adjoint unitaries, so $|X+Y|=I_{q^2}+XY$.
Since $\operatorname{tr}_q(uv)=0$, normalization of the second trace gives
\[
 (\operatorname{tr}_q\otimes\tau)(|X+Y|)=q,
 \qquad \|(u,v)\|_{S_1^q(\ell_2)}=q\sqrt2,
\]
where the second identity follows from Lemma~\ref{lem:Cnorm}.
This proves $A_1\geq\sqrt2$ whenever at least two variables are selected.
For one variable, $\|C\otimes\epsilon\|_1=\|C\|_1$ gives $A_1=1$.
Finally, if two sites are full, the six-variable witness in
Proposition~\ref{prop:quantum-rademacher-lower} gives $A_1\geq\sqrt{8/3}$.
\end{proof}
\end{example}
\subsection{Scalar coefficients and unbounded difference multiplicity}

\begin{example}[Scalar coefficients with $Z_2=\infty$]\label{ex:scalar-unbounded-differences}
The requirement that the inequality hold for all matrix coefficients
in Theorem~\ref{thm:quantitative-converse} is essential. In the scalar
notation fixed in the introduction, for $V\subset\Z$ we have
\[
 \alpha_1(V)
 =\sup_{0\ne(c_n)\in c_{00}(V)}
 \frac{\left(\sum_{n\in V}|c_n|^2\right)^{1/2}}
 {\displaystyle\int_0^1
       \left|\sum_{n\in V}c_ne^{2\pi int}\right|\,dt},
\]
where $c_{00}(V)$ denotes the finitely supported complex scalar
families. The paired lacunary sets considered below have infinitely
many representations of the difference $1$, but their scalar
constants are strictly less than $\sqrt2$.
\end{example}

\begin{lemma}[The independent paired Steinhaus model]\label{lem:paired-Steinhaus}
Let $S$ and $(S_k)$ be independent Steinhaus variables. For finitely
supported complex scalar families $(a_k)$ and $(b_k)$, put
$B=\sum_k(|a_k|^2+|b_k|^2)$. Then
\begin{equation}\label{eq:scalar-paired-independent}
 \mathbb E\left|\sum_k(a_k+b_kS)S_k\right|
 \geq\sqrt{\frac2\pi}\,B^{1/2}.
\end{equation}
\end{lemma}

\begin{proof}
We may assume $B>0$. Conditioning on $S$ and applying Sawa's
sharp Steinhaus inequality~\cite{Sawa1985} gives
\[
 \mathbb E\left|\sum_k(a_k+b_kS)S_k\right|
 \geq\frac{\sqrt\pi}{2}\,
       \mathbb E_S\left(\sum_k|a_k+b_kS|^2\right)^{1/2}.
\]
Put $T=\sum_k\overline{a_k}b_k$, so $|T|\leq B/2$ and the
expression under the square root is $B+2\operatorname{Re}(ST)$.
The function
\[
 \phi(r)=\int_0^{2\pi}\sqrt{1+r\cos\theta}\,
                    \frac{d\theta}{2\pi},\qquad -1\leq r\leq1,
\]
is even and concave. Its minimum is therefore
$\phi(1)=2\sqrt2/\pi$. Substitution proves
\eqref{eq:scalar-paired-independent}.
\end{proof}

\begin{lemma}[A dyadic coupling]\label{lem:dyadic-coupling}
Fix $m\geq3$, put $q_k=2^{k^2}$, and let $t$ be uniform on $[0,1)$.
Write $z=e^{2\pi it}$ and $\chi_k=z^{q_k}$ for $k\geq m$.
On an enlarged probability space there are independent Steinhaus
variables $S,(S_k)_{k\geq m}$ such that, for every finitely supported
pair of scalar families $(a_k),(b_k)$,
\begin{equation}\label{eq:dyadic-coupling-error}
 \mathbb E\left|
   \sum_{k\geq m}(a_k+b_kz)\chi_k
   -\sum_{k\geq m}(a_k+b_kS)S_k\right|
 \leq e_m\left(\sum_{k\geq m}(|a_k|^2+|b_k|^2)\right)^{1/2},
\end{equation}
where
\begin{equation}\label{eq:dyadic-tail-error}
 e_m=\frac1{\sqrt6}
      \left(2\pi2^{-m^2}+\frac{4\pi}{\sqrt{15}}4^{-m}\right).
\end{equation}
\end{lemma}

\begin{proof}
Adjoin independent uniform variables $U_0,(U_k)_{k\geq m}$ on
$[0,1)$, independent of $t$. Put $L_k=2^{2k+1}$ and define
\[
 S=\exp\left(2\pi i\frac{\lfloor2^{m^2}t\rfloor+U_0}{2^{m^2}}\right),
 \qquad
 S_k=\exp\left(2\pi i\frac{\lfloor L_k\{q_kt\}\rfloor+U_k}{L_k}\right),
\]
where $\{x\}$ is the fractional part of $x$.
Outside the null set of dyadic rationals, the retained binary digits
occupy positions $1,\ldots,m^2$ for $S$ and
$k^2+1,\ldots,(k+1)^2$ for $S_k$. These blocks are disjoint.
Together with the independent uniform tails, this proves the
asserted independence and distributions. Moreover, $S$ is independent
of the entire error vector $(\chi_k-S_k)_{k\geq m}$, which uses
only later digits and the auxiliary tails $(U_k)_{k\geq m}$.

The difference of two independent uniform variables on $[0,1)$
has second moment $1/6$. Since $u\mapsto e^{2\pi iu}$ is
$2\pi$-Lipschitz, the coupling gives
\[
 \|z-S\|_2\leq\frac{2\pi2^{-m^2}}{\sqrt6},\qquad
 \|\chi_k-S_k\|_2\leq\frac{\pi4^{-k}}{\sqrt6}.
\]
Set
\[
 F=\sum_{k\geq m}(a_k+b_kz)\chi_k,\qquad
 F_1=\sum_{k\geq m}(a_k+b_kS)\chi_k,\qquad
 F_2=\sum_{k\geq m}(a_k+b_kS)S_k,
\]
and write $B=\sum_{k\geq m}(|a_k|^2+|b_k|^2)$.
Cauchy--Schwarz and orthogonality of the characters $\chi_k$ give
\[
 \mathbb E|F-F_1|
 \leq\|z-S\|_2\left(\sum_{k\geq m}|b_k|^2\right)^{1/2}
 \leq\frac{2\pi2^{-m^2}}{\sqrt6}\,B^{1/2}.
\]
Pointwise Cauchy--Schwarz, independence of $S$ and the error vector,
and $\mathbb E S=0$ give
\[
 \begin{aligned}
 \mathbb E|F_1-F_2|
 &\leq\left(\sum_{k\geq m}\|\chi_k-S_k\|_2^2\right)^{1/2}
          \left(\mathbb E\sum_{k\geq m}|a_k+b_kS|^2\right)^{1/2}\\
 &\leq\frac{4\pi}{\sqrt{90}}\,4^{-m}B^{1/2}.
 \end{aligned}
\]
Adding the two estimates proves~\eqref{eq:dyadic-coupling-error}.
\end{proof}

\begin{proposition}[Sharp dyadic tails]\label{prop:scalar-dyadic-sharp}
For $m\geq3$ let
\[
 V_m=\{2^{k^2},2^{k^2}+1:k\geq m\},
\]
and let $e_m$ be as in~\eqref{eq:dyadic-tail-error}. Then
$Z_2(V_m)=\infty$ and
\begin{equation}\label{eq:scalar-sharp-tails}
 \sqrt{\frac\pi2}\leq \alpha_1(V_m)
 \leq\left(\sqrt{\frac2\pi}-e_m\right)^{-1}<\sqrt2.
\end{equation}
Moreover,
$\alpha_1(V_3)<1.296$,
$\alpha_1(V_4)<1.262$, and
$\lim_{m\to\infty}\alpha_1(V_m)=\sqrt{\pi/2}$.
\end{proposition}

\begin{proof}
The repeated difference $1=(2^{k^2}+1)-2^{k^2}$ gives
$Z_2(V_m)=\infty$. Every polynomial supported on $V_m$ has the
form $F=\sum_{k\geq m}(a_k+b_kz)\chi_k$, and its coefficient
$\ell_2$ norm is $B^{1/2}$, where
$B=\sum_{k\geq m}(|a_k|^2+|b_k|^2)$.
Lemmas~\ref{lem:paired-Steinhaus} and~\ref{lem:dyadic-coupling}
give $\mathbb E|F|\geq(\sqrt{2/\pi}-e_m)B^{1/2}$, proving the
upper bound. The sequence $e_m$ decreases, and substitution at
$m=3,4$ gives the two stated numerical bounds and the strict
inequality in~\eqref{eq:scalar-sharp-tails}.

For the reverse estimate, use the normalized polynomials
\[
 F_{M,n}(z)=\frac{1+z}{\sqrt{2n}}
               \sum_{k=M}^{M+n-1}z^{2^{k^2}},\qquad M\geq m.
\]
Their coefficients have $\ell_2$ norm one. By
Lemma~\ref{lem:dyadic-coupling}, the corresponding independent
model is
\[
 G_{M,n}=\frac{1+S}{\sqrt2}\,
             \frac1{\sqrt n}\sum_{k=M}^{M+n-1}S_k,
 \qquad \mathbb E|F_{M,n}-G_{M,n}|\leq e_M.
\]
The complex central limit theorem gives
$n^{-1/2}\sum_kS_k\Rightarrow g$, where $g$ is a circular complex
Gaussian with $\mathbb E|g|^2=1$ and
$\mathbb E|g|=\sqrt\pi/2$. The second moments of the normalized
sums equal one, so their absolute values are uniformly integrable.
Independence of $S$ and $(S_k)$, together with
$\mathbb E|1+S|=4/\pi$, therefore gives
\[
 \lim_{n\to\infty}\mathbb E|G_{M,n}|
 =\frac4{\pi\sqrt2}\,\frac{\sqrt\pi}{2}
 =\sqrt{\frac2\pi}.
\]
Since $F_{M,n}$ is supported on $V_m$, it follows that
\[
 \alpha_1(V_m)
 \geq\left(\sqrt{\frac2\pi}+e_M\right)^{-1}
 \qquad(M\geq m).
\]
Letting $M\to\infty$ proves the lower bound. The upper bound and
$e_m\to0$ give the asserted limit.
\end{proof}

\begin{corollary}[An explicit scalar counterexample]\label{prop:scalar-unbounded-differences}
Let $V=\{2^{k^2},2^{k^2}+1:k\geq4\}$. Then $Z_2(V)=\infty$, whereas
\begin{equation}\label{eq:scalar-unbounded-differences}
 \alpha_1(V)
 \leq\left(\sqrt{\frac2\pi}
       -\frac1{\sqrt6}\left(\frac\pi{32768}
                      +\frac\pi{64\sqrt{15}}\right)\right)^{-1}
 <1.262<\sqrt2.
\end{equation}
\end{corollary}

\begin{proof}
Apply Proposition~\ref{prop:scalar-dyadic-sharp} with $m=4$ and
note that
$e_4=(\pi/32768+\pi/(64\sqrt{15}))/\sqrt6$.
\end{proof}

\subsection{Further group examples}

\begin{example}[Free subsets]\label{ex:free-subsets}
If $V$ is a left or right translate of a free set of generators in a
discrete group, then $Z_2^c(V),Z_2^r(V)\leq1$. Thus
Corollary~\ref{thm:Main-Group} gives $A_1(V)\leq\sqrt2$.
For a singleton $V$, both $A_1(V)$ and its scalar counterpart
$\alpha_1(V)$ equal one.
\end{example}
\begin{example}(Planar circles in $\Z^2$).
For $k\in\N$, let
$V_k=\{(m,n)\in\Z^2:m^2+n^2=k\}$.
Then $Z_2(V_k)\leq2$, as observed by Zygmund~\cite{Zygmund1974}:
for a fixed nonzero difference, its representations correspond to the
intersection of a circle with a distinct translate of itself, which
contains at most two points.
\begin{corollary}
Fix $k\in\N$ and a family $(C_{m,n})_{(m,n)\in V_k}$ in $M_d(\C)$.
Then
\begin{equation}
 \|(C_{m,n})_{(m,n)\in V_k}\|_{S_1^d(\ell_2)}
 \leq\sqrt3\left\|
    \sum_{(m,n)\in V_k}C_{m,n}e^{2\pi i(mt_1+nt_2)}
          \right\|_{L_1(\mathbb T^2;S_1^d)}.
\end{equation}
\end{corollary}
\end{example}
\begin{remark}[The Furstenberg set]\label{rem:furstenberg}
The product set $V=\{2^m3^n:m,n\geq0\}$ has larger difference
multiplicities than the preceding examples. In fact, $Z_2(V)\geq7$,
since
\[
 \begin{aligned}
 27-3=32-8=36-12=48-24&=24,\\
 72-48=96-72=216-192&=24.
 \end{aligned}
\]
These are seven distinct pairs in $V\times V$. Multiplying every
entry by $6$ gives the same lower bound when both exponents are
required to be positive.
\end{remark}

\subsection{Gaussian systems}

\begin{example} (Complex Gaussian variables).
Let $W=(\gamma_j)_{j\geq1}$ be independent standard complex Gaussian
variables, normalized by $\mathbb E\gamma_j=0$ and
$\mathbb E|\gamma_j|^2=1$. Thus their real and imaginary parts are
independent centered real Gaussians of variance $1/2$.
Their centered products have $Z_2 (W)=1$.
 
Haagerup and Musat~\cite[Theorems 2.1 and 2.9]{Haagerup2007} already gave $A_1(W)=\sqrt2 $.
\end{example}
\begin{example} (Real Gaussian variables).
For independent standard real Gaussians $W=(g_j)_{j\geq1}$, the
centered-product Gram matrix has norm $Z_2^c(W)=Z_2^r(W)=2$.
For distinct $w,x\in W$, the two pairings $(y,z)=(w,x)$ and $(x,w)$
are responsible for this value. The sharp Rademacher result of
Section~\ref{rad26:sec:main} and the central limit theorem give $A_1(W)=\sqrt2$.
\begin{corollary}\label{cor:real-gaussian-sharp}
For independent standard real Gaussians $(g_j)$ and any
$C_1,\ldots,C_l\in M_d(\C)$,
\begin{equation}\label{eq:real-gaussian-sharp}
 \left\|\sum_{j=1}^l C_j\otimes g_j\right\|_{L_1(\Omega;S_1^d)}
 \leq\|(C_j)\|_{S_1^d(\ell_2)}
 \leq\sqrt2\left\|\sum_{j=1}^l C_j\otimes g_j
                      \right\|_{L_1(\Omega;S_1^d)}.
\end{equation}
The constant $\sqrt2$ is optimal uniformly over $d$ and $l$.
\end{corollary}
\begin{proof}
We first record an invariance of the coefficient norm. If
$a_1,\ldots,a_r\in\C$ satisfy $\sum_q|a_q|^2=1$, then
\begin{equation}\label{eq:coefficient-replication}
 \|(a_qC_j)_{j,q}\|_{S_1^d(\ell_2)}
 =\|(C_j)_j\|_{S_1^d(\ell_2)}.
\end{equation}
Indeed, replication sends each column or row decomposition of
$(C_j)$ to one with the same corresponding square-function norms.
Conversely, the map $(D_{j,q})\mapsto(\sum_q\overline{a_q}D_{j,q})_j$
is a contraction on both the column and row spaces by
Cauchy--Schwarz, and is a left inverse to replication. Taking the
infimum over decompositions proves the equality.

Let $(\varepsilon_{j,q})$ be independent Rademacher variables and put
\[
 g_j^{(m)}=m^{-1/2}\sum_{q=1}^m\varepsilon_{j,q},\qquad
 F_m=\sum_{j=1}^l C_j g_j^{(m)}.
\]
Apply Theorem~\ref{rad26:thm:main} to the coefficients
$(m^{-1/2}C_j)_{j,q}$ and use
\eqref{eq:coefficient-replication}. This gives
\[
 \mathbb E\|F_m\|_1\leq\|(C_j)\|_{S_1^d(\ell_2)}
 \leq\sqrt2\,\mathbb E\|F_m\|_1.
\]
The multivariate central limit theorem gives
$(g_1^{(m)},\ldots,g_l^{(m)})\Rightarrow(g_1,\ldots,g_l)$.
Moreover,
\[
 \mathbb E\|F_m\|_1^2
 \leq d\,\mathbb E\operatorname{tr}(F_m^*F_m)
 =d\sum_{j=1}^l\operatorname{tr}(C_j^*C_j).
\]
Thus $(\|F_m\|_1)_m$ is uniformly integrable. Passing to the limit
proves \eqref{eq:real-gaussian-sharp}.

Finally, $(g_{2j-1}+ig_{2j})/\sqrt2$ are independent standard
complex Gaussians. Apply any proposed real Gaussian inequality to
the paired coefficients $C_j/\sqrt2$ and $iC_j/\sqrt2$.
Their coefficient norm equals that of $(C_j)$ by
\eqref{eq:coefficient-replication}, so the same constant would hold
for complex Gaussians. Haagerup--Musat's sharpness result for complex
Gaussians~\cite[Theorems~2.1 and~2.9]{Haagerup2007} therefore forces
that constant to be at least $\sqrt2$.
\end{proof}
\end{example}

\section{An open question for subsets of the integers}\label{sec:integer-Sidon-question}
Every infinite subset $V\subset\Z$ satisfies $A_1(V)\geq\sqrt2$.
More generally, \eqref{eq:infinite-abelian-lower} gives the same lower
bound for infinite subsets of any discrete abelian group, by the
Haagerup--Musat sharpness argument recalled in
Lemma~\ref{lem:commuting-lower}.
Together with Theorem~\ref{thm:Main}, this shows that $Z_2(V)=1$
implies $A_1(V)=\sqrt2$ for infinite integer sets.
On the other hand, Theorem~\ref{thm:quantitative-converse} gives only
$Z_2(V)\leq2$ when $A_1(V)=\sqrt2$. This leaves the following question.

\begin{question}\label{question:integer-Sidon}
Let $V$ be an infinite subset of $\Z$. Is it true that
\begin{equation}\label{eq:integer-Sidon-question}
 A_1(V)=\sqrt2\quad\Longleftrightarrow\quad Z_2(V)=1?
\end{equation}
\end{question}

\begin{remark}[Sidon sets in the sense of Erd\H{o}s]\label{rem:Erdos-Sidon}
The condition $Z_2(V)=1$ is equivalent to $V$ being a Sidon set in the
sense of Erd\H{o}s, also called an additive $B_2$ set: every integer has
at most one representation $a+b$ with $a,b\in V$ and $a\leq b$.
Repetitions $a=b$ are allowed in this definition.
Equivalently, $a+b=c+d$ forces equality of the two unordered pairs,
counting multiplicities. This follows by rewriting a nontrivial
sum identity as two distinct representations of a nonzero difference,
and conversely. This additive notion should be distinguished from
the Sidon property in harmonic analysis.
The unresolved implication in \eqref{eq:integer-Sidon-question} is
$A_1(V)=\sqrt2\Rightarrow Z_2(V)=1$; by
Theorem~\ref{thm:quantitative-converse}, it suffices to settle the
case $Z_2(V)=2$.
\end{remark}

\begin{remark}[Why infinitude is necessary]\label{rem:finite-three-point}
The implication $A_1(V)=\sqrt2\Rightarrow Z_2(V)=1$ fails for finite
sets. Indeed, let $E=\{0,1,2\}$. The difference $1$ has the two
representations $1-0=2-1$, and no nonzero difference has more than
two representations, so $Z_2(E)=2$. Nevertheless,
\begin{equation}\label{eq:three-point-sharp}
 \|(C_0,C_1,C_2)\|_{S_1^d(\ell_2)}
 \leq\sqrt2\int_{\mathbb T}
       \operatorname{tr}|C_0+C_1z+C_2z^2|\,dm(z),
 \qquad A_1(E)=\sqrt2,
\end{equation}
where $m$ is normalized Haar measure. We include a verification of
both the upper bound and its sharpness.

For a triple $D=(D_0,D_1,D_2)$ in the unit ball of
$S_\infty^d(\ell_2)$, consider the block Hankel matrix
\[
 \mathsf H(D)=
 \begin{pmatrix}
 D_0&D_1&D_2\\
 D_1&D_2&0\\
 D_2&0&0
 \end{pmatrix}.
\]
We claim that $\|\mathsf H(D)\|\leq\sqrt2$.
First assume $D_k=D_k^*$, so that $\sum_{k=0}^2D_k^2\leq I_d$.
For $\xi=(x,y,z)\in(\C^d)^3$, put
\[
 u=(x,y),\quad v_0=(x/2,0),\quad
 v_1=(y/2,x/2),\quad v_2=(z,y/2),
 \qquad a=\|x\|^2+\|y\|^2,\quad b=\|z\|^2.
\]
Since $\sum_k\|(D_k\oplus D_k)u\|^2\leq a$ and
$\sum_k\|v_k\|^2=a/2+b$, Cauchy--Schwarz gives
\[
 \begin{aligned}
 |\langle\mathsf H(D)\xi,\xi\rangle|
 &=\left|2\operatorname{Re}\sum_{k=0}^2
          \langle(D_k\oplus D_k)u,v_k\rangle\right|\\
 &\leq2\sqrt{a(a/2+b)}
 \leq\sqrt2(a+b)=\sqrt2\|\xi\|^2.
 \end{aligned}
\]
This proves the claim in the self-adjoint case.
For general $D_k$, apply it to
$\widetilde D_k=\left(\begin{smallmatrix}0&D_k\\D_k^*&0\end{smallmatrix}\right)$.
The row and column bounds give $\sum_k\widetilde D_k^2\leq I_{2d}$,
and a permutation of coordinates identifies
$\mathsf H(\widetilde D)$ with
$\left(\begin{smallmatrix}0&\mathsf H(D)\\
\mathsf H(D)^*&0\end{smallmatrix}\right)$.
Thus the same norm bound holds for $\mathsf H(D)$.

Now put $F(z)=C_0+C_1z+C_2z^2$.
The standard matrix-valued analytic factorization
$H_1(S_1^d)=H_2(S_2^d)\cdot H_2(S_2^d)$, used in
\cite{LustPiquard1991}, gives, for every $\varepsilon>0$, analytic
factors $F=GH$ with
$\|G\|_{L_2(S_2^d)}\|H\|_{L_2(S_2^d)}
\leq\|F\|_{L_1(S_1^d)}+\varepsilon$.
Writing $G(z)=\sum_{j\geq0}G_jz^j$ and
$H(z)=\sum_{k\geq0}H_kz^k$, Hilbert--Schmidt duality and Parseval give
\[
 \begin{aligned}
 \left|\sum_{n=0}^2\operatorname{tr}(C_nD_n^*)\right|
 &=\left|\sum_{\substack{j,k\geq0\\j+k\leq2}}
       \operatorname{tr}(H_kD_{j+k}^*G_j)\right|\\
 &\leq\|\mathsf H(D^*)\|
       \left(\sum_{j=0}^2\|G_j\|_2^2\right)^{1/2}
       \left(\sum_{k=0}^2\|H_k\|_2^2\right)^{1/2}\\
 &\leq\sqrt2\bigl(\|F\|_{L_1(S_1^d)}+\varepsilon\bigr).
 \end{aligned}
\]
Here $D^*=(D_0^*,D_1^*,D_2^*)$.
Letting $\varepsilon\downarrow0$ and taking the supremum over $D$
in the dual formula \eqref{eq:coefficient-duality} proves the upper bound in
\eqref{eq:three-point-sharp}.

For sharpness, take
\[
 C_0=e_{11},\qquad C_1=e_{12}+e_{21},\qquad C_2=e_{22}.
\]
Both square-function sums equal $2I_2$, so
Lemma~\ref{lem:Cnorm} gives
$\|(C_0,C_1,C_2)\|_{S_1^2(\ell_2)}=2\sqrt2$.
On the other hand,
\[
 F(z)=\begin{pmatrix}1&z\\z&z^2\end{pmatrix}
     =\begin{pmatrix}1\\z\end{pmatrix}\begin{pmatrix}1&z\end{pmatrix}
\]
has trace norm $2$ for every $z\in\mathbb T$. Hence
$\|F\|_{L_1(S_1^2)}=2$, and the ratio is exactly $\sqrt2$.
Thus $Z_2(E)=2$ but $A_1(E)=\sqrt2$, showing why infinitude is
necessary in Question~\ref{question:integer-Sidon}.
\end{remark}

\appendix
\section{Proof of Lemma~\ref{rad26:lem:kkt}}\label{app:cadilhac-kkt}
The decomposition in Lemma~\ref{rad26:lem:kkt} is due to
Cadilhac~\cite[Theorem~4.2 and Remark~4.6]{Cadilhac19}.
We include the following finite-dimensional separation proof of
its formulation for a prescribed dual maximizing family.

\begin{proof}[Proof of Lemma~\ref{rad26:lem:kkt}]
Let $(C_k)$ and $(D_k)$ be as in the lemma, and put
$S=\sum_kD_k^2$. The proof is based on the Slater--KKT argument; see
\cite[Sections~5.9.1--5.9.2]{BoydVandenberghe2004}.
Let $\mathcal H_d$ be the real vector space of self-adjoint matrices
in $M_d(\C)$. Define, for $E_k \in \mathcal H_d$,
\[
 \ell(E)=\sum_k\operatorname{tr}(C_kE_k),\qquad
 \psi(E)=I_d-\sum_kE_k^2.
\]
The map $\psi$ is concave in the Loewner order. Indeed, for
$E,F\in\mathcal H_d^n$ and $0\leq\theta\leq1$,
\[
 \psi(\theta E+(1-\theta)F)
 -\theta\psi(E)-(1-\theta)\psi(F)
 =\theta(1-\theta)\sum_k(E_k-F_k)^2\geq0.
\]

Let
\[
 \mathcal U=\{(H,t)\in\mathcal H_d\times\R:
       H<\psi(E),\ t<\ell(E)\text{ for some }E\in\mathcal H_d^n\}.
\]
This set is nonempty and open, and the concavity of $\psi$ together
with the linearity of $\ell$ shows that it is convex.
Moreover, $(0,\ell(D))\notin\mathcal U$: otherwise some $E$ would satisfy
$\psi(E)>0$ and $\ell(E)>\ell(D)$, contradicting the maximality of $D$.
By the Hahn--Banach separation theorem, there is a nonzero element
$(R,\beta)\in\mathcal H_d\times\R$ such that
\[
 \operatorname{tr}(RH)+\beta t\leq\beta\ell(D)
 \qquad((H,t)\in\mathcal U).
\]
If $(H,t)\in\mathcal U$, then
$(H-rxx^*,t-s)\in\mathcal U$ for every column vector $x\in\C^d$ and all scalars $r,s\geq0$.
Letting $r$ and $s$ tend separately to infinity in the separating
inequality shows that $x^*Rx\geq0$ for every $x$ and that $\beta\geq0$.
Thus $R\geq0$.

For any $E=(E_k)\in\mathcal H_d^n$ and $\varepsilon>0$, insert
$H=\psi(E)-\varepsilon I_d$ and $t=\ell(E)-\varepsilon$ into the
separating inequality, and then let $\varepsilon\downarrow0$. We obtain
\begin{equation}\label{pineq}
 \operatorname{tr}(R\psi(E))+\beta\ell(E)\leq\beta\ell(D)
 \qquad(E\in\mathcal H_d^n).
\end{equation}
If $\beta=0$, taking $E=0$ and using $\psi(0)=I_d$ gives
$\operatorname{tr}(R)\leq0$. Since $R\geq0$, this implies $R=0$,
a contradiction. Hence $\beta>0$.

Set $P=R/\beta\geq0$ and define
\[
 L(E,P)=\ell(E)+\operatorname{tr}(P\psi(E)).
\]
Inequality~\eqref{pineq} says that $L(E,P)\leq\ell(D)$ for every
$E\in\mathcal H_d^n$. In particular, $L(D,P)\leq\ell(D)$.
Since $P\geq0$ and $\psi(D)\geq0$, equality holds, and therefore
\[
 L(E,P)\leq\ell(D)=L(D,P)\qquad(E\in\mathcal H_d^n).
\]
Consequently, $D$ maximizes $L(E,P)$ over the whole real vector
space $\mathcal H_d^n$,
and
\begin{equation}\label{rad26:eq:complementarity}
 \operatorname{tr}\bigl(P(I_d-S)\bigr)=0.
\end{equation}
Since $P\geq0$ and $S\leq I_d$, this equality gives
\[
 \|P^{1/2}(I_d-S)^{1/2}\|_2^2
 =\operatorname{tr}(P(I_d-S))=0.
\]
Multiplying the zero product by $P^{1/2}$ and $(I_d-S)^{1/2}$
shows that $P(I_d-S)=0$. Taking adjoints gives $(I_d-S)P=0$.
Thus $SP=PS=P$; no invertibility of $P$ is needed.

We now derive the first identity in \eqref{rad26:eq:kkt}. Fix an index $k$ and
an arbitrary self-adjoint matrix $H=H^*$. Let
\[
 E_k(s)=D_k+sH,\qquad E_j(s)=D_j\quad(j\neq k).
\]
The scalar function $g(s)=L(E(s),P)$ has a maximum at $s=0$, so
$g'(0)=0$. Since
\[
 (D_k+sH)^2=D_k^2+s(D_kH+HD_k)+s^2H^2,
\]
we obtain
\begin{align*}
 0=g'(0)
 &=\left.\frac{d}{ds}L(E(s),P)\right|_{s=0}\\
 &=\operatorname{tr}(C_kH)
   -\operatorname{tr}\bigl(P(D_kH+HD_k)\bigr)\\
 &=\operatorname{tr}\bigl((C_k-PD_k-D_kP)H\bigr).
\end{align*}
By the arbitrarity of $H$, we get
\[
 C_k=PD_k+D_kP\qquad(1\leq k\leq n).
\]

Finally, if $P=0$, then $C_k=PD_k+D_kP=0$ for every $k$. Hence
$P\neq0$ whenever $(C_k)\neq0$.
\end{proof}

\section*{Acknowledgment}
The authors are grateful to \'{E}ric Ricard for insightful suggestions in the study of the sharp
Rademacher constant.  They also thank \'{E}ric Ricard, L\'{e}onard Cadilhac,
Javier Parcet, Yin Sheng and Ping Zhong for helpful discussions.
OpenAI's ChatGPT assisted in editing and in determining constants
for examples.
The third author is supported by NSF grant DMS-2625657.

\clearpage
\begingroup
\small
\let\originalthebibliography\thebibliography
\renewcommand{\thebibliography}[1]{%
  \originalthebibliography{#1}%
  \setlength{\itemsep}{0pt plus .2pt}%
  \setlength{\parskip}{0pt}}


\endgroup
\end{document}